\documentclass[twoside,leqno,10pt]{article}
\usepackage{graphics,ifthen}
\usepackage{latexsym}
\usepackage{amsmath}
\usepackage{amssymb}
\usepackage{mathrsfs}

\begin{document}
\input{latexP.sty}
\input{referencesP.sty}
\input epsf.sty

\def\ind{\stackrel{\mathrm{ind}}{\sim}}
\def\iid{\stackrel{\mathrm{iid}}{\sim}}

\def\Definition{\stepcounter{definitionN}\
    \Demo{Definition\hskip\smallindent\thedefinitionN}}
\def\EndDefinition{\EndDemo}
\def\Example#1{\Demo{Example [{\rm #1}]}}
\def\EndExample{\qed\EndDemo}
\def\Category#1{\centerline{\Heading #1}\rm}
\
\def\e{\text{\hskip1.5pt e}}
\newcommand{\eps}{\epsilon}
\newcommand{\proof}{\noindent {\bf Proof:\ }}
\newcommand{\remarks}{\noindent {\bf Remarks:\ }}
\newcommand{\note}{\noindent {\bf Note:\ }}
\newcommand{\examp}{\noindent {\bf Example:\ }}
\newcommand{\Lower}[2]{\smash{\lower #1 \hbox{#2}}}
\newcommand{\ben}{\begin{enumerate}}
\newcommand{\een}{\end{enumerate}}
\newcommand{\bi}{\begin{itemize}}
\newcommand{\ei}{\end{itemize}}
\newcommand{\hp}{\hspace{.2in}}

\newtheorem{lw}{Proposition 3.1, Lo and Weng (1989)}
\newtheorem{thm}{Theorem}[section]
\newtheorem{defin}{Definition}[section]
\newtheorem{prop}{Proposition}[section]
\newtheorem{lem}{Lemma}[section]
\newtheorem{cor}{Corollary}[section]
\newcommand{\rb}[1]{\raisebox{1.5ex}[0pt]{#1}}
\newcommand{\mc}{\multicolumn}
\newcommand{\Acr}{\mathscr{A}}
\newcommand{\Bcr}{\mathscr{B}}
\newcommand{\Ucr}{\mathscr{U}}
\newcommand{\Gcr}{\mathscr{G}}
\newcommand{\Dcr}{\mathscr{D}}
\newcommand{\CS}{\mathscr{C}}
\newcommand{\Fcr}{\mathscr{F}}
\newcommand{\Icr}{\mathscr{I}}
\newcommand{\Lcr}{\mathscr{L}}
\newcommand{\Mcr}{\mathscr{M}}
\newcommand{\Ncr}{\mathscr{N}}
\newcommand{\Pcr}{\mathscr{P}}
\newcommand{\Qcr}{\mathscr{Q}}
\newcommand{\Scr}{\mathscr{S}}
\newcommand{\Tcr}{\mathscr{T}}
\newcommand{\Xcr}{\mathscr{X}}
\newcommand{\Vcr}{\mathscr{V}}
\newcommand{\Ycr}{\mathscr{Y}}
\newcommand{\qcr}{\mathscr{q}}
\newcommand{\scr}{\mathscr{s}}
\newcommand{\indic}{\mathbb{I}}
\newcommand{\E}{\mathbb{E}}
\newcommand{\F}{\mathbb{F}}
\newcommand{\I}{\mathbb{I}}
\newcommand{\Q}{\mathbb{Q}}
\newcommand{\X}{\mathbb{X}}
\newcommand{\Pe}{\mathbb{P}}
\newcommand{\M}{\mathbb{M}}
\newcommand{\Wbb}{\mathbb{W}}

\def\Beta{\text{Beta}}
\def\Dir{\text{Dirichlet}}
\def\DP{\text{DP}}
\def\P{{\bf p}}
\def\fhat{\widehat{f}}
\def\GA{\text{gamma}}
\def\ind{\stackrel{\mathrm{ind}}{\sim}}
\def\iid{\stackrel{\mathrm{iid}}{\sim}}
\def\J{{\bf J}}
\def\K{{\bf K}}
\def\min{\text{min}}
\def\N{\text{N}}
\def\p{{\bf p}}
\def\U{{\bf U}}
\def\W{{\bf W}}
\def\S{{\bf S}}
\def\T{{\bf T}}
\def\y{{\bf y}}
\def\t{{\bf t}}
\def\m{{\bf m}}
\def\X{{\bf X}}
\def\Y{{\bf Y}}
\def\tps{\mbox{\scriptsize ${\theta H}$}}   
\def\ups{\mbox{\scriptsize ${P_{\theta}(g)}$}}   
\def\vps{\mbox{\scriptsize ${\theta}$}}   
\def\vups{\mbox{\scriptsize ${\theta >0}$}}   
\def\hps{\mbox{\scriptsize ${H}$}}   
\def\rps{\mbox{\scriptsize ${(\theta+1/2,\theta+1/2)}$}}   
\def\sps{\mbox{\scriptsize ${(1/2,1/2)}$}}   

\newcommand{\reals}{{\rm I\!R}}
\newcommand{\PR}{{\rm I\!P}}
\def\Z{{\bf Z}}
\def\yy{{\mathcal Y}}
\def\rr{{\mathcal R}}
\def\BP{\text{beta}}
\def\ts{\tilde{t}}
\def\js{\tilde{J}}
\def\gs{\tilde{g}}
\def\fs{\tilde{f}}
\def\ys{\tilde{Y}}
\def\ps{\tilde{\mathcal {P}}}

\def\Report{Lancelot F. James}
\def\Author{Gamma Tilting and Bessel Occupation Times}
\pagestyle{myheadings}
\markboth{\Author}{\Report}
\thispagestyle{empty}

\bct\Heading  Gamma Tilting Calculus for GGC and Dirichlet means
with applications to Linnik processes and Occupation Time Laws for
Randomly Skewed Bessel Processes and Bridges. \lbk\lbk\smc
Lancelot F. James\footnote{\eightit \rm Supported in
part by grants HIA05/06.BM03 and DAG04/05.BM56 of the HKSAR.
\eightit Corresponding authors address.
                \rm The Hong Kong University of Science and Technology,
Department of Information and Systems Management, Clear Water Bay,
Kowloon, Hong Kong.
\rm lancelot\at ust.hk\\
\indent\eightit Keywords and phrases.
                \rm
          Bessel Process,
          Brownian Motion,
          Dirichlet Process,
          Generalized gamma convolution,
          L\'evy process,
          two-parameter Poisson-Dirichlet process
          }
\lbk\lbk \BigSlant The Hong Kong University of Science and
Technology\rm \lbk 
\ect \Quote This paper explores various interfaces between the
class of generalized gamma convolution~(GGC) random variables,
Dirichlet process mean functionals and phenomena connected to the
local time and occupation time of p-skew Bessel processes and
bridges discussed in Barlow, Pitman and Yor~(1989), Pitman and Yor
(1992, 1997). First, some general calculus for GGC and Dirichlet
process means functionals is developed. It then proceeds, via an
investigation of positive Linnik random variables, and more
generally random variables derived from compositions of  a stable
subordinator with GGC subordinators, to establish various
distributional equivalences between these models and phenomena
connected to local times and occupation times of what are defined
as randomly skewed Bessel processes and bridges. This yields a
host of interesting identities and explicit density formula for
these models. Randomly skewed Bessel processes and bridges may be
seen as a randomization of their p-skewed counterparts, and are
shown to naturally arise via exponential tilting. As a special
result it is shown that the occupation time of a p-skewed random
Bessel process or (generalized) bridge is equivalent in
distribution to the occupation time of a non-trivial randomly
skewed process.\EndQuote
\rm
\section{Introduction}
This paper explores various interfaces between the class of
generalized gamma convolution~(GGC) random variables, Dirichlet
process mean functionals and phenomena connected to the local time
and occupation time of skew Bessel processes and bridges discussed
in Barlow, Pitman and Yor~(1989), Pitman and Yor (1992, 1997b).
The latter concepts are also related to the now ubiquitous
two-parameter Poisson-Dirichlet models, see for instance Pitman
and Yor~(1997a) and Pitman~(2006). Specifically, we take the
viewpoint as in James~(2006), that GGC random variables can be
represented as linear functionals of Gamma processes and hence
there is a strong link to the distributional theory of Dirichlet
mean functionals developed in Cifarelli and Regazzini~(1990).
Here, we first develop two types of distributional results, one is
a simple but quite useful result involving scaling general
Dirichlet mean functionals by beta random variables, and the other
is a form of gamma tilting calculus. We then investigate the class
of positive Linnik random variables(see Devroye (1996, 1990)),
establishing new distributional identities, and  a link to local
and occupation times of skew Bessel processes and bridges. As a
nautral extension, we investigate
 compositions of a stable subordinator and a subordinator which
 has a GGC law. We use this to develop additionally a
 scaling calculus for (generalized) Dirichlet mean functionals. These
 points then lead to  the development of a series of interesting results, including explicit density formula, for what we call
randomly skewed Bessel processes and bridges. As a special case,
we establish equivalences in law between occupation times of
$p$-skewed Bessel processes and bridges and non-trivial randomly
skewed processes.  Some recent work related to our line of
investigation include, James~(2006), James, Lijoi and
Pr\"unster~(2006), Bertoin, Fujita, Roynette and Yor~(2006) and
Fujita and Yor~(2006). With some further reference to these latter
two works, we note that a key result in this manuscript is to show
that the random variable $X_{\alpha}=S_{\alpha}/S'_{\alpha}$,
where $S_{\alpha}$ and $S'_{\alpha}$ are iid $\alpha$-stable
random variables is indeed a Dirichlet mean functional. The
explicit distribution of $X_{\alpha}$ was obtained in
Lamperti~(1958) and, as shall be shown, this fact combined with
the correspondence to Dirichlet mean functionals allows one to
easily obtain explicit densities for a host of  models connected
to Bessel phenomena. For some other references relevant to skew
Bessel processes see Bertoin and Yor~(1996), Kasahara and
Watanabe~(2005) and Watanabe~(1995).

\section{The relationship between GGC, Gamma processes and Dirichlet mean functionals}
A positive infinitely divisible random variable $Z_{\theta}$ is
said to be a generalized Gamma convolution~(GGC) if for positive
$\lambda$ one can write its Laplace transform in the form
$$
\E[{\mbox e}^{-\lambda Z_{\theta}}]={\mbox e}^{-\theta
\vartheta(\lambda)}
$$
where for a sigma-finite measure $\nu$ on $(0,\infty),$
$$
\vartheta(\lambda)=\int_{0}^{\infty}\log(1+\lambda
x)\nu(dx)=\int_{0}^{\infty}\int_{0}^{\infty}(1-{\mbox
e}^{-s\lambda})s^{-1}{\mbox e}^{-s/x}\nu(dx).
$$
$\nu$ is chosen so that $\vartheta(\lambda)<\infty,$ and is often
referred to as a Thorin measure. For precise conditions  on $\nu$
see Bondesson~(1992). The corresponding L\'evy measure of
$Z_{\theta}$ is,
$$
\theta s^{-1}\int_{0}^{\infty}{\mbox e}^{-s/x}\nu(dx).
$$
We will say $Z_{\theta}$ is GGC$(\theta,\nu).$ Naturally every GGC
can be used to generate a subordinator which we denote as
$Z_{\theta}(t)$ for each $t\ge 0.$ We will also be interested in
the special case where $\nu=H$ is a probability measure, in that
case we write GGC$(\theta, H)$ and will sometimes refer to
$Z_{\theta}$ as an FGGC or finite GGC, a term we used in
James~(2006). As in that work, we will exploit a close connection
between GGC, and Gamma processes, and in the case of FGGC, its
connection with the Dirichlet process.[See Ferguson~(1973) for the
Dirichlet process and its use in Bayesian nonparametric
statistics]. Now let $\Gamma_{\theta \nu}$ denote a gamma process
on $(0,\infty)$ with shape parameter $(\theta \nu),$ this is a
completely random measure whose law is characterized by its
Laplace functional for some positive function $g$ as
$$
\E[{\mbox e}^{-\lambda\Gamma_{\theta \nu}(g)}]={\mbox e}^{-\theta
\int_{0}^{\infty}\log(1+\lambda g(x))\nu(dx)},
$$
where $\Gamma_{\theta \nu}(g):=\int_{0}^{\infty}g(x)\Gamma_{\theta
\nu}(dx).$ Hence setting $g(x)=x$, we see that every
GGC$(\theta,\nu)$ random variable can be written as a
mean~(linear) functional of $\Gamma_{\theta \nu}.$ Let
$(J_{k,\theta})$ denote the ranked jump sizes of a subordinator
based on a \textsc{Gamma}$(\theta)$ distribution, such that
$\sum_{k=1}^{\infty}J_{k,\theta}=G_{\theta}.$ Now recall that for
each $\theta>0$, $G_{\theta},$ which is a \textsc{Gamma}$(\theta)$
random variable, is independent of the sequence of ranked
probabilities $({\tilde P}_{k,\theta}=J_{k,\theta}/G_{\theta}).$
The law of the sequence $({\tilde P}_{k,\theta})$ is sometimes
referred to as the Poisson-Dirichlet distribution and, using
notation from Pitman and Yor~(1997a), is denoted PD$(0,\theta).$
These points imply that
$$
\int_{0}^{\infty}x\Gamma_{\theta
\nu}(dx)\overset{d}=G_{\theta}M_{\theta}(\nu)
$$
where $M_{\theta}(\nu)\overset{d}=\sum_{k=1}^{\infty}{\tilde
P}_{k,\theta}W_{k}$, where the $(W_{k})$ are independent of
$(J_{k,\theta})$ and are the points of a Poisson random measure
with mean intensity $\nu.$ Furthermore, the distribution of
$M_{\theta}(\nu)$ is characterized by its Cauchy-Stieltjes
transform of order $\theta$, which is
$$
\E[{(1+\lambda M_{\theta}(\nu))}^{-\theta}]=\E[{\mbox e}^{-\lambda
G_{\theta}M_{\theta}(\nu)}]={\mbox e}^{-\theta\vartheta(\lambda)}.
$$
Hence, importantly, it follows that
$Z_{\theta}\overset{d}=G_{\theta}M_{\theta}(\nu).$ Now when
$\nu=H$ we may define a Dirichlet process with shape parameter
$\theta H$ as $ P_{0,\theta}(\cdot)=\sum_{k=1}^{\infty}{\tilde
 P}_{k,\theta}\delta_{Z_{k}}(\cdot), $ where the $(Z_{k})$ are iid
with common distribution $H.$ Then a Dirichlet mean functional is
defined as
$$M_{\theta}(H)\overset{d}=\int_{0}^{\infty}xP_{0,\theta}(dx).$$
Hence, if $Z_{\theta}$ is a GGC$(\theta,H)$ then
$Z_{\theta}\overset{d}=G_{\theta}M_{\theta}(H).$

\Remark Note of course that when $\nu$ has infinite mass the use
of $\theta$ is somewhat redundant. However, it will play an
important role as the decomposition $G_{\theta}M_{\theta}(\nu)$
always makes sense.\EndRemark

\Remark In Bondesson~(1992) the L\'evy exponent of a GGC is
represented in a slightly different, but certainly equivalent, way
as
$$
\vartheta(\lambda)=\int_{0}^{\infty}\log(1+\lambda/y)U(dy).
$$
where $U$ is some sigma-finite measure. In this work, we are using
a representation which more closely matches that used in the
literature on Dirichlet process mean functionals. \EndRemark

\Remark See Vershik, Yor and Tsilevich~(2001) and James~(2005) for
some uses of the decomposition $G_{\theta}M_{\theta}(H)$\EndRemark

\Remark Throughout, for $0<\alpha<1,$ we will define $S_{\alpha}$
to be a unilateral $(\alpha)-$ stable random variable having the
Laplace transform
$$
\E[{\mbox e}^{-\lambda S_{\alpha}}]={\mbox e}^{-\lambda^{\alpha}}.
$$
This random variable is well known to be a GGC with an infinite
Thorin measure. See example 3.2.1 of Bondesson~(1992). \EndRemark
\section{Some calculus for GGC/Dirichlet mean functionals}
Suppose that $X$ has distribution $H$, and define the function
$$
\Phi(t)=\int_{0}^{\infty}\log(|t-x|)I(t\neq
x)H(dx)=\E[\log(|t-X|)I(t\neq X)]
$$
furthermore, letting $H(t)=\int_{0}^{t}H(dx)$, define
$$
\Delta_{\theta}(t|H)=\frac{1}{\pi}\sin(\pi \theta H(t)){\mbox
e}^{-\theta \Phi(t)}.
$$
Cifarelli and Regazzini~(1990, 1994), see also Cifarelli and
Mellili~(2000), apply inversion formula to obtain the
distributional formula for $M_{\theta}(H)$ as follows. For all
$\theta>0$, the cdf can be expressed as
\begin{equation}
\label{DPcdf} \int_{0}^{x}{(x-t)}^{\theta-1} \Delta_{\theta
}(t|H)dt
\end{equation}
provided that $\theta H$ possesses no jumps of size greater than
or equal to one.  If we let $f_{M_{\theta}}(\cdot|H)$ denote the
density of $M_{\theta}(H),$ it takes it simplest form for
$\theta=1$, which is \Eq \label{M1}
f_{M_{1}}(x|H)=\Delta_{1}(x|H)=\frac{1}{\pi}\sin(\pi H(x)){\mbox
e}^{-\Phi(x)}. \EndEq Formula for $\theta>1$ are also available
from Cifarelli and Regazzini~(1990,1994) and Cifarelli and
Mellili~(2000). However, we shall not explicitly use this and
instead give an expression for the density which holds for all
$\theta>0$, recently obtained by James, Lijoi and
Pr\"unster~(2006), as follows,

\Eq \label{generaldensity} f_{M_{\theta}}(x|H)=
 \int_{0}^{x}{(x-t)}^{\theta-1} d_{\theta}(t|H)dt
\EndEq where $$d_{\theta}(t|H)=\frac{d}{dt}\sin(\pi \theta
H(t)){\mbox e}^{-\theta \Phi(t)}.$$

We do point out, that except for some of the very recent work we
mention here, there are very few examples of mean functionals
where the explicit density has been calculated. This is due in
part to the fact that it is not necessarily obvious how to write
$\Delta_{\theta }(t|H)$ in a nice form, or otherwise how to
calculate $\Phi(x).$ Our results in this section provide some
additional tools to obtain more explicit expressions and also
point to apparently unknown inter-relationships between different
mean functionals which occurs through a tilting operation.

The first result, in Theorem 3.1 below,  shows that multiplying a
mean functional by a certain beta random variable may lead to more
simplified expressions. Additionally, we  point out that this
result was constructed in an effort to better understand the
result in Proposition 15 of Pitman and Yor~(1997b) and, as we
show, paves the way for other interesting results related to
occupation times and local times of Bessel processes.

\Remark First, we introduce a bit more notation. The notation
$G_{a}$, for $a>0,$ will denote a Gamma random variable with shape
$a$ and scale $1,$ with law denoted as \textsc{Gamma}(a).
Throughout we will use the notation $B_{a,b}$ to denote a Beta
random variable with parameters $(a,b)$, and write
\textsc{Beta}$(a,b)$ to denote this law. Additionally, unless
otherwise specified, if for random variables $X$ and $Y$ we write
the product XY, it will be assumed that $X$ and $Y$ are
independent. Additionally the notation $X'$ means that $X'$ is an
independent random variable with the same distribution as $X$
\EndRemark

\begin{thm}Let $X$ denote a random variable with distribution $H.$
For $0<p\leq 1$, let $Y_{p}$ denote a \textsc{Bernoulli}$(p)$
random variable independent of $X.$ Then the distribution of
$XY_{p}$ is given by $H^{(p)}(dv):=p H(dv)+(1-p)\delta_{0}(dv).$
Now define $\Phi(t):=\E[\log(|t-X|)\indic(X\neq t)].$ Hence
$\Phi^{(p)}(t):=\E[\log(|t-XY_{p}|)\indic(XY_{p}\neq
t)]=p\Phi(t)+(1-p)\log(t)\indic(t\neq 0).$ Then for $\theta>0,$
\Eq \label{betaid} M_{\theta}(H^{(p)})\overset{d}=M_{p \theta
}(H)B_{\theta p,\theta(1-p)} \EndEq and
$G_{\theta}M_{\theta}(H^{(p)})\overset{d}=G_{\theta p}M_{\theta
p}(H).$ Hence GGC$(\theta, H^{(p)})$=GGC$(\theta p,H)$ for all
$\theta>0$and $0<p<1.$ As special cases,
\begin{enumerate}
\item[(i)]$M_{1}(H^{(p)})=M_{p}(H)B_{p,1-p},$ with density, for
$0<p<1,$ \Eq \label{densimple}\frac{1}{\pi}x^{-(1-p)}\sin(\pi
p[1-H(x)]){\mbox e}^{-p\Phi(x)} \EndEq
\item[(ii)]$M_{\frac{1}{p}}(H^{(p)})=M_{1}(H)B_{1,\frac{(1-p)}{p}}.$
\end{enumerate}
\end{thm}
\Proof First note the distributional identity, $G_{\theta
p}\overset{d}=G_{\theta}B_{\theta p,\theta(1-p)}$, for all
$\theta>0$. Then for~\mref{betaid} it suffices to check the
equivalence $G_{\theta}M_{\theta}(H^{(p)})\overset{d}=G_{\theta
p}M_{\theta p}(H).$ But this is obvious upon taking Laplace
transforms. The density in~\mref{densimple} follows as a special
case of~\mref{M1}. Note additionally the identity
$\sin(\pi[pH(x)+1-p])=\sin(\pi[p(1-H(x)]),$ for $0<p<1.$\EndProof
\subsubsection{Remark about obtaining the distribution of a FGGC subordinator
for all time points} It is well known that if $\{Z(t),t\ge 0\}$ is
some subordinator, even if we have an explicit expression for its
density or distribution at time say 1, i.e. the random variable
$Z(1),$ it is not obvious how to find explicit densities or
distribution functions for general $Z(t),$ at fixed time points
$t.$ Of course, by the convolution properties of infinitely
divisible random variables one can describe the law of a random
variable for $Z(t)$ for $t>1,$ in terms of a linear combination of
independent random variables based on $\{Z(s):s\leq 1\}.$ Hence,
from this point of view,  one may simply concentrate on obtaining
explicit expressions for the density of random variables for
$Z(t)$, where $t$ is some fixed value in $(0,1].$ In the case
where $\{Z_{\theta}(t);t\ge 0\}$ is a FGGC subordinator, one may
obtain an explicit description of its density for each fixed
$0<t\leq 1$, provided that $\Phi$ has a tractable form, using
Theorem 3.1. Formally if $0<\theta t<1$, it follows that
$$
Z_{\theta}(t)\overset{d}=G_{\theta t}M_{\theta
t}(H)\overset{d}=G_{1}B_{\theta t,1-\theta t}M_{\theta
t}(H)\overset{d}=G_{1}M_{1}(H^{(\theta t)})
$$
where $M_{1}(H^{(\theta t)})$ has a density (of non-integral
form),
$$
\frac{1}{\pi}x^{\theta t-1}\sin(\pi \theta t[1-H(x)]){\mbox
e}^{-\theta t\Phi(x)}.
$$
The point here is that density of $M_{1}(H^{(\theta t)})$
generally has a simpler description than that of $M_{\theta
t}(H).$ Naturally, these models agree when $\theta t=1.$ \Remark
These points and Theorem 3.1 are connected to an expression for a
GGC$(\beta,H)$ density for $0<\beta<1$ given in Bondesson~(1992,
p. 37-38). Relationships to beta random variables are not noted
there. \EndRemark
\subsection{Some basic exponential tilting calculus for Gamma mixtures}
Here we develop distributional results for exponential tilting of
scale mixtures of Gamma random variables. We then specify this to
the case of exponential tilting of GGC models which always can be
written in this form. An important point is that we are able to
establish a clear distributional link between a $M_{\theta}(\nu)$
and a corresponding random variable obtained from tilting
$G_{\theta}M_{\theta}(\nu),$ which, as a by-product, has new
implications for the study of Dirichlet mean functionals These
results will play a fundamental role throughout the text.
\begin{thm} Suppose that $W=G_{\theta}M$ where $M$ and $G_{\theta}$ are
independent positive random variables and $G_{\theta}$ is
\textsc{Gamma}$(\theta)$. Then for $b$ and $c$ positive numbers,
the random variable, ${\tilde W}_{\theta}$ with density ${\mbox
e}^{-c/b w}b^{-1}f_{W}(w/b)/E[{\mbox e}^{-cG_{\theta}M}]$
satisfies,
$$
{\tilde W}_{\theta}\overset{d}=(b/c)G_{\theta}{Y}_{\theta,c}
$$
where,
\begin{enumerate}
\item[(i)]${Y}_{\theta,c}$ is a random variable with density
$$
\frac{f_{Y}(y){[(1-y)]}^{\theta}}{\E[{(1+cM)}^{-\theta}]},$$ where
$Y\overset{d}=cM/(cM+1),$ and $\E[{(1+cM)}^{-\theta}]=\E[{\mbox
e}^{-cG_{\theta}M}].$
\item[(ii)]Equivalently the distribution of $Y_{\theta,c}$ is the
distribution of $cM/(cM+1)$ taken with respect to the density
$(cm+1)^{-\theta}f_{M}(m)/\E[{\mbox e}^{-cG_{\theta}M}].$
\item[(iii)]Conversely
$f_{M}(x)={(1+x)}^{\theta-2}f_{Y_{\theta,1}}(\frac{x}{1+x})\E[{(1+M)}^{-\theta}]$
\end{enumerate}
\end{thm}
\Proof Let us proceed by checking Laplace and Cauchy transforms.
By direct argument, working with the density, the Laplace
transform of ${\tilde W}_{\theta}$ evaluated at $\lambda$ is,
$$
\frac{\E[{(1+(c+\lambda b)M)}^{-\theta}]}{\E[{(1+cM)}^{-\theta}]}
$$
which, for $c\neq 0,$ may be re-written as
$$
\frac{\E[{(1+c(1+\lambda
(b/c))M)}^{-\theta}]}{\E[{(1+cM)}^{-\theta}]}
$$
Now for simplicity we set $b/c=1$ and note that the Laplace
transform of $G_{\theta}Y_{\theta,c}$ is, by definition of $Y_{
\theta,c},$
$$
\E[(1+\lambda
Y_{\theta,c})^{-\theta}]=\int_{0}^{\infty}{\left(1+\lambda
\frac{cx}{cx+1}\right)}^{-\theta}\frac{{(1+cx)}^{-\theta}f_{M}(x)}{\E[{(1+cM)}^{-
\theta}]}dx
$$
Simple algebra completes the result. \EndProof The next two
results are specialized to GGC.
\begin{prop}Let $L$ be a GGC$(\theta,\nu)$ random variable with density denoted $f_{L}.$ Hence
$L=G_{\theta}M_{\theta}(\nu)$. Now define a tilted random variable
$\tilde{L}$ having density proportional to ${\mbox e}^{-c/b
w}b^{-1}f_{W}(w/b)$ for $c\ge 0$ and $b>0.$ Denote this law of
${\tilde L}$ as GGC$^{(b,c)}(\theta,\nu).$ It then follows that
for $c>0$ ${\tilde L}\overset{d}=(b/c)L_{c}$ where $L_{c}$ is
GGC$^{(c,c)}(\theta,\nu).$
\begin{enumerate}
\item[(i)]The L\'evy measure of ${\tilde L}$ is
$$
s^{-1}{\mbox e}^{-(c/b)s}\int_{0}^{\infty}e^{-s/br}\nu(dr)
$$
\item[(ii)]When $c=0$, $\tilde{L}=bL$. Otherwise for $c>0$,
${\tilde L}\overset{d}=(b/c)G_{\theta}M_{\theta}(\nu^{(c,c)})$.
The Thorin measure $\nu^{(c,c)}$ is determined by $\nu$ as
indicated by the following expression for the L\'evy density of
$G_{\theta}M_{\theta}(\nu^{(c,c)}),$
$$
s^{-1}\int_{0}^{\infty}e^{-s\frac{cr+1}{cr}}\nu(dr)=s^{-1}\int_{0}^{\infty}e^{-s/y}\nu^{(c,c)}(dy)
$$
\item[(iii)]Referring to Theorem 3.2,
$M_{\theta}(\nu^{(c,c)})\overset{d}=Y_{\theta,c}$, for
$M\overset{d}=M_{\theta}(\nu).$
\end{enumerate}
\end{prop}
\Proof The first statement is just a standard result for
exponentially tilting infinitely divisible random variables
coupled with the specific form of a GGC L\'evy measure. Statement
[(ii)] is just algebra, which yields importantly the
representation $G_{\theta}M_{\theta}(\nu^{(c,c)}).$ Statement
[(iii)] then follows from Theorem 3.2\EndProof

 We close this section with an  important variation in the case where
$\nu$ is a probability measure.
\begin{prop}Let $L$ be a GGC$(\theta, H)$, then there exists a
random variable $X$ with distribution $H$ and a random variable
$A^{*}_{c}=cX/cX+1$ on $[0,1]$ such that if $L_{c}$ is now
GGC$^{(c,c)}(\theta,H)$ it is equivalently GGC$(\theta,Q_{c})$,
where $Q_{c}$ is the distribution of $A^{*}_{c}$ derived from H
and is a special case of $\nu^{(c,c)}.$ In other words the L\'evy
measure of $L_{c}$ may be written as
$$
\theta s^{-1}\E[{\mbox e}^{-s/A^{*}_{c}}]
$$
 Additionally
the following distributional relationships hold.
\begin{enumerate}
\item[(i)]Suppose that the density of $M_{\theta}(H)$, say $f_{M_{\theta}}(\cdot|H)$
is known. Then the density of $M_{\theta}(Q_c)$ is expressible as
$$
f_{M_{\theta}}(y|Q_{c})=\frac{{(1-y)}^{\theta-2}}{c\E[{(1+cM_{\theta}(H))}^{-\theta}]}f_{M_{\theta}}\left(\frac{y}{c(1-y)}|H\right)$$
\item[(ii)]Conversely, if the density of $M_{\theta}(Q_{c})$,
$f_{M_{\theta}}(\cdot|Q_{c}),$ is known then the density of
$M_{\theta}(H)$ is given by
$$
f_{M_{\theta}}(x|H)={(1+x)}^{\theta-2}f_{M_{\theta}}\left(\frac{x}{1+x}|Q_{1}\right)\E[{(1+M_{\theta}(H))}^{-\theta}]$$
\end{enumerate}
\end{prop}
\Remark The gamma exponential tilting operation is quite special
and we point out that for a fixed value of $\theta$,
$Y_{\theta,c}$ cannot achieve all possible distributions on
$[0,1]$. So for example when $\theta=1$, $Y_{1,1}$ cannot be
\textsc{Uniform}$[0,1]$. This is equivalent to noting that a
density on $[0,1]$ of the form
$$
f_{Y}(y)\propto\frac{1}{(1-y)}
$$
does not exist. In general, the \textsc{Uniform}$[0,1]$ is
possible for $Y_{\theta,c}$, only when $0<\theta<1.$ \EndRemark
\subsection{An example with some connections to the occupation time of skew Brownian bridge}
Set $H=U[0,1]$ to denote that the corresponding random variable
$X$ is \textsc{Uniform}$[0,1].$ It is known from Diaconis and
Kemperman~(1994) that the density of $M_{1}(U[0,1])$ is \Eq
\label{DPUni} \frac{e}{\pi}\sin(\pi y)y^{-y}{(1-y)}^{-(1-y)}{\mbox
{ for }}0<y<1. \EndEq Note furthermore that $W=G_{1}M_{1}(U[0,1])$
is GGC$(1,U[0,1])$ and has a rather strange Laplace transform,
$$
\E[{\mbox e}^{-\lambda G_{1}M_{1}(U[0,1])}]={\mbox
e}{(1+\lambda)}^{-(\frac{\lambda+1}{\lambda})}.
$$

We can use this fact combined with the previous results to obtain
a new explicit expression for the density of what we believe
should be an important mean functional and corresponding
infinitely divisible random variable.
\begin{prop}Let $G_{1}/E$ be the ratio of two independent exponential $(1)$ random variables having density
$\zeta(dx)/dx=(1+x)^{-2}{\mbox { for }}x>0$. Now let $L$ denote a
GGC$(1,\zeta)$ random variable, with log Laplace transform
$\log\E[{\mbox e}^{-\lambda
G_{1}M_{1}(\zeta)}]=-\frac{\lambda}{\lambda-1}\log(\lambda).$ Then
equivalently the L\'evy measure of $L$ is given by
$$
s^{-1}\E[{\mbox e}^{-s\frac{E}{G_{1}}}]=s^{-1}\E[{\mbox
e}^{-s\frac{G_{1}}{E}}]
$$
Furthermore $L\overset{d}=G_{1}M_{1}(\zeta)$, where $M_{1}(\zeta)$
has density,
$$
f_{M_{1}}(x|\zeta)=\frac{1}{\pi}\sin(\pi
\frac{x}{1+x})x^{-\frac{x}{(1+x)}}{\mbox { for }}x>0.
$$
\end{prop}
\Proof First note that it is straightforward to show that
$\E[{\mbox e}^{-L}]=\E[{(1+M_{1}(\zeta))}^{-1}]=e^{-1}.$ This fact
also establishes the existence of $L.$ Now we see that
$G_{1}/(G_{1}+E)\overset{d}=U_{[0,1]}$. The result then follows by
applying statement [(ii)] of Proposition 3.2 to
\mref{DPUni}\EndProof

In view of the remarks in Section 3.0.1, we now give a description
of the laws of the subordinators associated with the two random
variables above.
\begin{prop}Suppose the $Z(t)$ is a subordinator where $Z(1)$ is
GGC$(1,U[0,1])$, then for $0<t<1$,
$Z(t)\overset{d}=G_{1}M_{1}(U^{(t)}[0,1]),$ where
$M_{1}(U^{(t)}[0,1])$ has density,
$$
\frac{e^{t}}{\pi}\sin(\pi t(1-
y))y^{t(1-y)-1}{(1-y)}^{-t(1-y)}{\mbox { for }}0<y<1.
$$
If $Z(t)$ is such that $Z(1)$ is GGC$(1,\zeta)$ then then for
$0<t<1$, $Z(t)\overset{d}=G_{1}M_{1}(\zeta^{(t)}),$ where
$M_{1}(\zeta^{(t)})$ has density,
$$
\frac{1}{\pi}\sin(\pi \frac{t}{1+x})x^{\frac{t}{(1+x)}-1}{\mbox {
for }}x>0.
$$
\end{prop}
\Proof This now follows from Theorem 3.1, Proposition 3.3 and
~\mref{DPUni}. \EndProof

 Once we have the density in $M_{1}(\zeta)$ we can then
extend the result for $M_{1}(U[0,1])$ to that of
$M_{1}(O^{br}_{1/2,p}),$ where $O^{br}_{1/2,p}$ denotes the
distribution of the random variable
$$
A^{br}_{1/2,p}=\int_{0}^{1}\indic(B^{(1/2,1/2)}_{p}(s)>0)ds\overset{d}=\frac{p^{2}G_{1}}{p^{2}G_{1}+q^{2}E}.
$$
In the notation above $B^{(1/2,1/2)}_{p}(s)$ denotes a $p$-skew
Brownian motion. Hence $A^{br}_{1/2,p}$ denotes the time spent
positive by this process up till time $1.$ The following result is
otherwise not obvious.
\begin{prop}Define $A^{br}_{1/2,p}\overset{d}=p^{2}G_{1}/[p^{2}G_{1}+q
^{2}E].$ The random variable $A^{br}_{1/2,p}$ is equivalent in
distribution to the time spent positive of a $p$-skew Brownian
bridge having density
$O^{br}_{1/2,p}(dy)/dy=p^{2}q^{2}/[p^{2}(1-y)+q^{2}y].$ Now let
$L_{p}$ denote a GGC$(1,O^{br}_{1/2,p})$ random variable. Then
equivalently the L\'evy measure of $L_{p}$ is given by
$$
s^{-1}\E[{\mbox e}^{-s/A^{br}_{1/2,p}}]
$$
Furthermore $L_{p}\overset{d}=G_{1}M_{1}(O^{br}_{1/2,p})$, where
$M_{1}(O^{br}_{1/2,p})$ has density,
$$
f_{M_{1}(O^{br}_{1/2,p})}(y)=\frac{\kappa_{p}}{\pi}\sin\left(\frac{\pi
q^{2}y}{p^{2}(1-y)+q^{2}y}\right)y^{-\frac{
q^{2}y}{p^{2}(1-y)+q^{2}y}}{(1-y)}^{-\frac{
p^{2}(1-y)}{p^{2}(1-y)+q^{2}y}} {\mbox { for }}0<y<1.
$$
Where for $c=p^{2}/q^{2}$
$$
1/\kappa_{p}=c\E[(1+cM({\zeta})^{-1}]=c{\mbox
e}^{-\int_{0}^{\infty}\frac{\log(1+cx)}{{(1+x)}^{2}}dx}
$$

\end{prop}
\Remark The specific densities in Proposition 3.3 and 3.5 yield
the not immediately obvious identities.
$$
\Phi_{1}(x)=-\int_{0}^{\infty}\frac{\log(|x-y|)}{(1+y)^{2}}dy=-\frac{x}{1+x}\log(x)
$$
and, more so,
$$
\Phi_{2}(x)=-\int_{0}^{1}\frac{\log(|x-y|)p^{2}q^{2}}{[p^{2}(1-y)+q^{2}y]}dy=\log\left(\kappa_{p}y^{-\frac{
q^{2}y}{p^{2}(1-y)+q^{2}y}}{(1-y)}^{-\frac{
p^{2}(1-y)}{p^{2}(1-y)+q^{2}y}}\right)
$$
where these are appropriate versions of $\Phi.$ \EndRemark

\subsection{Reconciling some results of Cifarelli and Mellili}To
further illustrate our point we show how to reconcile two
apparently unrelated results given in Cifarelli and
Mellilli~(2000).  Let $\Lambda_{1/2,1/2}$ denote the distribution
of the arcsine law, that is a $B_{1/2,1/2}$ random variable.
Cifarelli and Melilli~(2000, p.1394-195) show that for all
$\theta>0$
$M_{\theta}(\Lambda_{1/2,1/2})\overset{d}=B_{\theta+1/2,\theta+1/2}.$
Now define the probability density
$$
\varrho_{1/2}(x)=\frac{1}{\pi}x^{-1/2}{(1+x)}^{-1}.
$$
Cifarelli and Mellili then show that for $\theta\ge 1$
$M_{\theta}(\varrho_{1/2})$ has the density proportional to
$$
x^{\theta-1/2}{(1+x)}^{-(\theta+1)}.
$$
Hjort and Ongaro~(2005) recently extend this result for all
$\theta>0$ and also note the normalizing constant appearing in
Cifarelli and Mellili~(2001) is incorrect. Here however we note
that if $X$ has density $\varrho_{1/2}$ then
$X\overset{d}=G_{1/2}/G'_{1/2}$ where $G_{1/2}$ and $G'_{1/2}$ are
independent and identically distributed gamma random variables.
Now using the known fact that
$$
B_{1/2,1/2}\overset{d}=\frac{G_{1/2}/G'_{1/2}}{G_{1/2}/G'_{1/2}+1},
$$
we see that $B_{1/2,1/2}$ is a special case of $A_{1}$ in
Proposition 3.2. It is now evident that one could use the result
$M_{\theta}(\Lambda_{1/2,1/2})\overset{d}=B_{\theta+1/2,\theta+1/2},$,
coupled with statement [(ii)] of Proposition 3.2, to easily obtain
the density of $M_{\theta}(\varrho_{1/2})$ for all $\theta>0$.
Similar to section 3.2 one could then use the density of
$M_{\theta}(\varrho_{1/2})$ to obtain results for mean functionals
based on the law of $p^{2}G_{1/2}/[p^{2}G_{1/2}+q^{2}G'_{1/2}].$
We will encounter this class of models again in the next coming
sections and see how they arise as laws of occupations times of
Brownian bridge and related models.

\section{A tilted positive Linnik process}
In this section we present details of a  subordinator, $Z$, which
is a FGGC such that $Z(1)$ has the Laplace transform, \Eq
\label{LaplaceLinnik}
{\left[\frac{1}{1+c^{\alpha}}(1+{(c+b\lambda)}^{\alpha})\right]}^{-\theta}
\EndEq for parameters $0<\alpha<1$, $\theta>0$, $c\ge 0$ and
$b>0.$ The first important thing to note is that $Z(t)$ is
equivalent in distribution to that of $Z(1)$ with the parameter
$\theta$ replaced by $\theta t.$ This as we shall show more
specifically, yields the desirable property of being able to
explicitly identify the distribution of $Z(t)$ for all $t.$ In
this generality we believe that this process has not been studied
in any detail, but setting $c=0,$ we see that the Laplace
transform becomes $$ {(1+b^{\alpha}\lambda^{\alpha})}^{-\theta}.$$

It is easy to verify that this corresponds to the random variable
defined as
$$bL_{\alpha,\theta}:=bG^{1/\alpha}_{\theta}S_{\alpha}.$$
The random variable, $L_{\alpha,\theta}$, is known in the
literature[see for instance Bondesson~(1992, p.38) and
Devroye~(1990, 1996)] and is sometimes called a positive Linnik
process and has a host of interesting properties and
distributional representations. Now let us introduce the random
variable which will play a key role throughout the remainder of
the paper.  Let $X_{\alpha}\overset{d}=S_{\alpha}/S'_{\alpha}$
denote the random variable having density \Eq
\varrho_{\alpha}(y)=\frac{\sin(\pi
\alpha)}{\pi}\frac{y^{\alpha-1}}{y^{2\alpha}+2y^{\alpha}\cos(\pi
\alpha)+1}{\mbox { for }}y>0 \label{denX}.\EndEq Then the L\'evy
exponent associated with $L_{\alpha,\theta}$ is
$$
\tilde {\psi}_{\alpha,\theta}(\omega):=
\theta\ln(1+\omega^{\alpha})=\int_{0}^{\infty}(1-{\mbox
 e}^{-\lambda s})l_{\theta,\alpha}(s)ds
$$
where $$l_{\theta,\alpha}(s)=\frac{\alpha\theta }{s}\phi_{\alpha}(
s^{\alpha})= \frac{\alpha \theta}{s} \E[{\mbox
e}^{-sX_{\alpha}}]=\alpha\theta s^{-1}\E[{\mbox
e}^{-s/X_{\alpha}}]$$ is the L\'evy density of the Linnik process.
Specifically, $$ \phi_{\alpha}(q)=\E[{\mbox
e}^{-qS^{-\alpha}_{\alpha}}]=\sum_{k=0}^{\infty}\frac{1}{\Gamma(1+k\alpha)}{(-q)}^{k}=\E[{\mbox
e}^{-q^{1/\alpha}X_{\alpha}}]
$$
is the Mittag Leffler function, which equates with the Laplace
transform of $S^{-\alpha}_{\alpha}.$ The equivalences involving
$X_{\alpha}$ are probably not that well known but serve to
identify the density of $X_{\alpha}$~\mref{denX} as the Thorin
measure, $\nu$ of $L_{\alpha,\theta}.$ This Thorin measure is
identified by Bondesson~((1992, p. 38) and confirms the fact that
$L_{\alpha,\theta}$ is an FGGC. That is,
$$
\tilde {\psi}_{\alpha,\theta}(\omega):=\theta \alpha
\int_{0}^{\infty}(1-{\mbox
 e}^{-\omega s})s^{-1}\int_{0}^{\infty}{\mbox
 e}^{-s/x}\varrho_{\alpha}(x)dxds.
$$
One of the things will be looking for are alternative
representations of the distribution of $L_{\alpha,\theta}.$ For
instance it is known from  Devroye~(1996), that when $\theta=1$
one has $L_{\alpha,1}=G_{1}X_{\alpha}$, where $G_{1}$ is
exponential $(1)$ and $X_{\alpha}$ has density ~\mref{denX}.
Summarizing, we see that $L_{\alpha,\theta}$ is a GGC$(\alpha
\theta,\varrho_{\alpha})$ and~\mref{LaplaceLinnik} corresponds to
the class GGC$^{(b,c)}(\alpha\theta,\varrho_{\alpha})$ We will
establish various distributional equivalences which then lead to
results for Dirichlet mean functionals and Bessel occupation
times. First we describe some more pertinent features of
$X_{\alpha}.$
\begin{prop}Let $X_{\alpha}\overset{d}=S_{\alpha}/S'_{\alpha},$ having
density~\mref{denX}. Then,
\begin{enumerate}
\item[(i)]The cdf of $X_{\alpha}$ can be represented explicitly as \Eq
\label{cdfXY}
F_{X_{\alpha}}(x)=1-\frac{1}{\pi\alpha}\cot^{-1}\left(\cot(\pi
\alpha)+\frac{x^{\alpha}}{\sin(\pi \alpha)}\right) \EndEq
\item[(ii)]Its inverse is given by \Eq \label{invX}
F^{-1}_{X_{\alpha}}(y)={\left[\frac{\sin(\pi \alpha(y))}{\sin(\pi
\alpha(1-y))}\right]}^{1/\alpha} \EndEq
\item[(iii)]The equations~\mref{cdfXY} and~\mref{invX} yield the quite useful
identity, \Eq \label{sinXid}
 \sin(\pi
\alpha(1-F_{X_{\alpha}}(y)))=y^{-\alpha}\sin(\pi \alpha
F_{X_{\alpha}}(y))=\frac{\sin(\pi\alpha)}{{[y^{2\alpha}+2y^{\alpha}\cos(\pi
\alpha)+1]}^{1/2}} \EndEq
\end{enumerate}
\end{prop}
\Proof This derivation of the cdf is influenced by arguments in
Fujita and Yor~(2006) where it becomes clear that it is easier to
work with the density of $[X_{\alpha}]^{1/\alpha}.$ The cdf of
this quantity is easily obtained. Statements [(ii)] and [(iii)]
then become apparent. \EndProof

\Remark There are several things interesting about the alternative
representation of the L\'evy measure involving $X_{\alpha}$.
First, from Chaumont and Yor~(2003, sec. 4.19 and 4.21),. the
distributional equivalence can also be seen to arise from the fact
that $G^{1/\alpha}_{1}\overset{d}=G_{1}/S'_{\alpha}$, where
$S'_{\alpha}$ is a stable random variable independent of
$S_{\alpha}$, and the fact that
$X_{\alpha}\overset{d}=S_{\alpha}/S'_{\alpha}$. The density of
$X_{\alpha}$ can be traced back to the work of Lamperti~(1958) but
arises later elsewhere. One notes also that the density for
$X_{\alpha}$ has a simple form as compared with a stable law of
index $0<\alpha<1$.  \EndRemark \Remark One notes that the random
variable $L_{\alpha, \theta}$ is conditionally a stable random
variable and hence is heavy-tailed. So this limits its practical
applicability to a variety of problems. A natural way to create a
random variable which has moments is to exponential tilt the
density of $L_{\alpha,\theta}.$ This is precisely how the Laplace
transform, ~\mref{LaplaceLinnik}, comes about. It is easy to see
that,
$$
\E[{\mbox e}^{-\lambda Z(1)}]=\frac{\E[{\mbox e}^{-(b\lambda+c)
L_{\alpha,\theta}}]}{\E[{\mbox e}^{-c L_{\alpha,\theta}}]}
$$
In other words the density of $Z(1)$ is given by
$$
f_{Z(1)}(x)=\frac{1}{b}(1+c^{\alpha}){\mbox
e}^{-xc/b}f_{L_{\alpha,\theta}}(x/b)=\frac{{\mbox e}^{-\frac
{x}{b}(p/q)^{1/\alpha}}}{b(1-p)}f_{L_{\alpha,\theta}}(x/b).
$$
\EndRemark
\subsection{Connection to occupation times of Bessel bridges and $P_{\alpha,\theta}(C)$}
Let $J_{1}\ge J_{2}\ge \cdots$ denote the ranked jump sizes of an
$(\alpha)$-stable subordinator such that
$\sum_{k=1}^{\infty}J_{k}=T_{\alpha,0}$ is a unilateral
$(\alpha)$-stable random variable. Recall from Pitman and
Yor~(1997a) that the sequence of ranked probabilities
$(P_{k}=J_{k}/T_{\alpha,0})$ is said to have a two-parameter
Poisson Dirichlet with specification $PD(\alpha,0).$ Furthermore,
for $\alpha\theta>0,$ one obtains the PD$(\alpha,\alpha \theta)$
law of $(P_{k})$ by mixing over the \emph{conditional} law of
$(P_{k})|T_{\alpha,0}=t$ with respect to the distribution of the
random variable, say $T_{\alpha,\alpha\theta}$, which has density
$\Pe(T_{\alpha,\alpha,\theta}\in dt)\propto
t^{-\alpha\theta}\Pe(T_{\alpha,0}\in dt).$  One may then
introduce, independent of $(P_{k}),$ a sequence of iid random
variables $(Z_{k})$ having some common non-atomic law. Then as in
Pitman~(1996), see also Ishwaran and James~(2001), analogous to
the Dirichlet process, one may define the class of
PD$(\alpha,\alpha\theta)$ random probability measures as
$P_{\alpha,\alpha\theta}(\cdot)=\sum_{k=1}^{\infty}P_{k}\delta_{{Z}_{k}}(\cdot).$

Now it will be shown that the exponential tilting operation to
obtain ~\mref{LaplaceLinnik} reveals a strong, albeit initially
unexpected, connection to work of Pitman and Yor~(1997b) and
Barlow, Pitman and Yor~(1989),  on occupation time models for skew
Bessel bridges and more generally laws of
$P_{\alpha,\alpha\theta}(C)$ for some set $C$ such that
$\E[P_{\alpha,\theta}(C)]=p.$

This is seen for the case where $c\neq 0,$ which allows one to to
rewrite ~\mref{LaplaceLinnik} as,
 \Eq
 \label{CSBesBridge}
\frac{1}{{(q+p{(1+\frac{b}{c}\lambda)}^{\alpha})}^{\theta}} \EndEq
for $p=c^{\alpha}/1+c^{\alpha}$ and $q=1-p.$ We recognize that for
$b=c$, this equates with the Cauchy-Stieljtes transform of order
$\alpha \theta$ of a two parameter $(\alpha,\alpha \theta)$
Poisson Dirichlet random probability measure evaluated at some set
$C$, denoted as $P_{\alpha,\alpha\theta}(C)$, where
$\E[P_{\alpha,\alpha\theta}(C)]=p.$ That is,
$$
\E[{(1+\lambda P_{\alpha,\alpha \theta}(C))}^{-\theta
\alpha}]=\frac{1}{{(q+p{(1+\lambda)}^{\alpha})}^{\theta}}.
$$

Furthermore, when $\theta=1$, this is the Cauchy-Stieltjes
transform of order $\alpha$ of the time spent positive up to time
$1$ of a $p$-skew Bessel bridge of dimension $2-2\alpha,$ as can
be seen from Pitman and Yor~(1997b, eq(75)). This random variable
can be represented as,
$$
A^{br}_{\alpha,p}\overset{d}=P_{\alpha,\alpha}(C)=\int_{0}^{1}\indic(B^{(\alpha,\alpha)}_{p}(s)>0)
ds
$$
where $B^{(\alpha,\alpha)}_{p}$ is the corresponding $p$-skew
Bessel bridge. We also note another important connection to
Lamperti~(1958) and Barlow, Pitman and Yor~(1989), Pitman and
Yor~(1997b). Noting that $c={(p/q)}^{1/\alpha}$, the random
variable \Eq \label{repA}
A_{\alpha,p}\overset{d}=P_{\alpha,0}(C)\overset{d}
=\frac{cX_{\alpha}}{1+cX_{\alpha}}\overset{d}=\frac{p^{1/\alpha}S'_{\alpha}}{p^{1/\alpha}S'_{\alpha}+q^{1/\alpha}S_{\alpha}}
\overset{d}=\int_{0}^{1} \indic(B^{(\alpha)}_{p}(s)>0)ds \EndEq
equates with the time spent positive by $B^{(\alpha)}_{p}$, now a
$p$-skew Bessel process of dimension $2-2\alpha$ with skewness
parameter $p$, up to time $1$. When $p=1/2$, one obtains the usual
Bessel processes. Lamperti~(1958) shows that the density of
$A_{\alpha,p}$ is \Eq \label{Denlamperti}
    \Lambda_{\alpha,p}(dx)/dx=\frac{p\,q\,\sin(\alpha\pi)\:x^{\alpha-1}\,
    (1-x)^{\alpha-1}}{\pi\:[q^2\, x^{2\alpha}+p^2(1-x)^{2\alpha}+2pq\:x^\alpha\,(1-x)^\alpha\,
    \cos(\alpha\pi)]},
    \EndEq
We see that setting $p=1/2,\alpha=1/2$, yields L\'evy's(1939)
famous result that the time spent positive by Brownian motion up
to time $1$ has the Arcsine distribution. That is, $A_{1/2,1/2}$
is \textsc{Beta}$(1/2,1/2).$ Now interestingly from ~\mref{cdfXY}
and ~\mref{invX} we obtain a closed form expression for the cdf
and quantile function of $A_{\alpha,p},$ as
 \Eq \label{cdfA}
F_{A_{\alpha,p}}(y)=1-\frac{1}{\pi\alpha}\cot^{-1}\left(\cot(\pi
\alpha)+\frac{qy^{\alpha}}{p{(1-y)}^{\alpha}\sin(\pi
\alpha)}\right) \EndEq and its inverse given by \Eq \label{invA}
F^{-1}_{A_{\alpha,p}}(y)=\frac{{\left[\frac{p\sin(\pi
\alpha(y))}{q\sin(\pi \alpha(1-y))}\right]}^{1/\alpha}}{1+
{\left[\frac{p\sin(\pi \alpha(y))}{q\sin(\pi
\alpha(1-y))}\right]}^{1/\alpha}}=\frac{F^{-1}_{X_{\alpha}}(y)}{1+F^{-1}_{X_{\alpha}(y)}}\EndEq

\Remark Hereafter, for $c\neq 0,$ we shall set $c^{\alpha}=p/q.$
Note that if $B(s)$ denotes in a generic sense a Bessel process or
bridge, then the interpretation of  $p$ is that the $p$-skewed
version of $B(s)$ has the property that
$$
p=\Pe(B(s)>0).
$$
\EndRemark

 Now we relate the more general class of
PD$(\alpha,\alpha\theta)$ models to occupation laws of processes
and their accompanying local times. Noting Pitman and Yor~(1992,
p. 332) let $(\ell^{(\alpha)}_{t},t\ge 0)$ denote the right
continuous local time of a Bessel process, and let
$S_{\alpha}(t)\overset{d}=S_{\alpha}t^{1/\alpha}$ denote an
$(\alpha)$-stable subordinator, which satisfies the identities in
law, \Eq \label{scale}
S_{\alpha}(s)=\inf\{t:\ell^{(\alpha)}_{t}>s\} {\mbox { and }},
\frac{\ell^{(\alpha)}_{t}}{t^{\alpha}}\overset{d}=\frac{s}{{(S_{\alpha}(s))}^{\alpha}}\overset{d}=\frac{1}{(S_{\alpha})^{\alpha}}.
\EndEq For our purposes, we may set $S_{\alpha}(1)=T_{\alpha,0},$
and furthermore
${[\ell^{(\alpha)}_{1}]}^{1/\alpha}\overset{d}=1/T_{\alpha,0}.$
Now using the scaling property~\mref{scale} one may construct
local times $\{\ell^{(\alpha,\alpha\theta)}_{t}:t\ge 0\}$  with
laws specified by
$\Pe({[\ell^{(\alpha,\alpha\theta)}_{1}]}^{1/\alpha}\in
ds)=\Pe(1/T_{\alpha,\alpha\theta}\in ds).$ As a special case,
$\ell^{(\alpha,\alpha)}$ denotes the local time of a Bessel
bridge, say $B^{(\alpha,\alpha)}.$ Associated with
$\{\ell^{(\alpha,\alpha\theta)}_{t}:t\ge 0\}$ are what we shall
call , \emph{generalized Bessel bridges} say
$B^{(\alpha,\alpha\theta)}(t)$ with law specified by
$\int_{0}^{\infty}\Pe(B^{(\alpha)}(t)|\ell^{(\alpha)}_{1}=s)\Pe(\ell^{(\alpha,\alpha\theta)}_{1}\in
ds).$ Such processes may be found in Definition 3.14 of Perman,
Pitman and Yor~(1992). Now letting $B^{(\alpha,\alpha\theta)}_{p}$
denote a $p$-skewed version of such processes we define their
times spent positive up to time $1$ as,
$A^{(\alpha,\alpha\theta)}_{\alpha,p}$ satisfying, $$
A^{(\alpha,\alpha\theta)}_{\alpha,p}\overset{d}=P_{\alpha,\alpha\theta}(C)
$$ Note that
$B^{(\alpha,\alpha\theta)}_{1/2}:=B^{(\alpha,\alpha\theta)}.$

 \Remark The representation of the
occupation time of a skew Bessel process in~\mref{repA} was given
by Barlow, Pitman and Yor~(1989). It will play a fundamental role
in our understanding of the construction of randomized versions
and related matters. In effect this boils down to the
interpretation of $c$ as $c^{\alpha}=p/q$ \EndRemark \Remark In
this remark we demonstrate how the positive Linnik may be
interpreted as the distribution of a time changed occupation time.
It is known, see Barlow, Pitman and Yor~(1989) or section 4 of
Pitman and Yor~(1992), that for $S_{\alpha}(t)$ an inverse local
time,
$$
A_{\alpha,p}(S_{\alpha}(t))=\int_{0}^{S_{\alpha}(t)}\indic(B^{(\alpha)}_{s,p}>0)
ds
\overset{d}=S_{\alpha}(p)t^{1/\alpha}\overset{d}=S_{\alpha}(pt)\overset{d}=p^{1/\alpha}t^{1/\alpha}S_{\alpha}
$$
is an $(\alpha)$-stable subordinator. Hence letting
$(\Gamma_{\theta}(t),t\ge 0)$ denote an independent gamma
subordinator satisfying for each fixed $t$,
$\Gamma_{\theta}(t)\overset{d}=G_{\theta t},$ it follows that
$$
A_{\alpha,p}(S_{\alpha}(\Gamma_{\theta}(t)))\overset{d}=S_{\alpha}(p\Gamma_{\theta}(t))\overset{d}=p^{1/\alpha}G^{1/\alpha}_{\theta
t}S_{\alpha}=p^{1/\alpha}S_{\alpha}(\Gamma_{\theta}(t)).
$$\EndRemark
\Remark
 The random variables
$T_{\alpha,\alpha\theta}^{-\alpha}$ are the $\alpha$-diversity of
the PD$(\alpha,\alpha \theta)$ exchangeable partitions as
described in Pitman~(2003, Proposition 13). Furthermore it may be
read from Perman, Pitman and Yor~(1992, p. 31) and Pitman~(2006)
that for $\theta\ge 1$ $T_{\alpha,\theta
\alpha}\overset{d}=T_{\alpha,\alpha\theta-\alpha}U_{\theta
\alpha,1-\alpha}$, where $U_{\theta \alpha,1-\alpha}$ is
independent of $T_{\alpha,\theta \alpha}$ but not of
$T_{\alpha,\alpha\theta-\alpha},$ and marginally $U_{\theta
\alpha,1-\alpha}\overset{d}=B_{\alpha\theta,1-\alpha}.$ We will
sometimes refer to all the relevant local times and occupations
times, already defined, as PD$(\alpha,\alpha\theta)$ processes.
\EndRemark
\subsection{Distributional results for the Linnik class and local times}
We first establish various distributional identities related to
the positive Linnik random variable. Importantly, we will show
that random variables based on the PD$(\alpha,\alpha\theta)$
models are Dirichlet mean functionals. \Remark Throughout we will
be using the fact that if $X$ is a gamma or $(\alpha)$-stable
random variable, then the independent random variables $X,Y,Z$
satisfying $XY\overset{d}=XZ$ imply that  $Y\overset{d}=Z.$ For
precise conditions see Chaumont and Yor~(2003, sec. 1.12 and
1.13). \EndRemark

\begin{thm} Let $L_{\alpha,\theta}=G^{1/\alpha}_{\theta}S_{\alpha}$ denote a generalized Linnik random variable i.e. a GGC$(\alpha\theta, \varrho_{\alpha})$.
Then, we have the distributional equivalences

\begin{itemize}\item[(i)] For all $\theta>0,$
$ L_{\alpha,\theta}\overset{d}=G_{\theta
\alpha}S_{\alpha}/{T_{\alpha,\alpha\theta}}\overset{d}=G_{\theta
\alpha}M_{\alpha \theta}(\varrho_{\alpha}),$ which implies
$$S_{\alpha}/T_{\alpha,\alpha\theta}\overset{d}=M_{\alpha
\theta}(\varrho_{\alpha})\overset{d}={[\ell^{(\alpha,\alpha\theta)}_{S_{\alpha}}]}^{1/\alpha}$$

\item[(ii)]In particular when $\theta=1$,
$
L_{\alpha,1}\overset{d}=G_{1}X_{\alpha}\overset{d}=G_{1}B_{\alpha,1-\alpha}M_{\alpha}(\varrho_{\alpha}).$
Hence
$$X_{\alpha}\overset{d}=B_{\alpha,1-\alpha}M_{\alpha}(\varrho_{\alpha})\overset{d}=M_{1}(\varrho^{(\alpha)}_{\alpha})$$
\item[(iii)]For $0<\theta<1,$
$$
L_{\alpha,\theta}\overset{d}=B^{1/\alpha}_{\theta,1-
\theta}G_{1}X_{\alpha}\overset{d}=G_{1}B_{\theta \alpha,1-\theta
\alpha}\frac{S_{\alpha}}{T_{\alpha,\alpha\theta}}\overset{d}=G_{\theta
\alpha+1-\alpha}B_{\alpha\theta,1-\alpha}\frac{S_{\alpha}}{T_{\alpha,\alpha\theta}}
$$
\item[(iv)]Statement (i) implies $G^{1/\alpha}_{\theta}\overset{d}=G_{\alpha
\theta}/T_{\alpha,\alpha\theta},$ for $\theta>0$. For
$0<\theta<1$,
$$
G^{1/\alpha}_{\theta}\overset{d}=\frac{G_{\alpha\theta}}{T_{\alpha,\alpha\theta}}\overset{d}=\frac{G_{1}B^{1/\alpha}_{\theta,1-\theta}}{S_{\alpha}}\overset{d}=\frac{G_{\alpha\theta+1-\alpha}}{T_{\alpha,\alpha\theta}}B_{\alpha\theta,1-\alpha}
$$
When $\theta=1$, then
$$
G^{1/\alpha}_{1}\overset{d}=\frac{G_{\alpha}}{S_{\alpha,\alpha}}\overset{d}=\frac{G_{1}}{S_{\alpha}}
$$
\end{itemize}
\end{thm}
\Proof The first equivalence in statement (i) is perhaps the most
non-obvious. The result is obtained by manipulating the density
representation of $G^{1/\alpha}_{\theta}S_{\alpha}$ as follows.
Let $f_{\alpha}$ denote the density of an $(\alpha)$-stable random
variable. Now note that the density of $L_{\alpha,\theta}$ is
obviously expressible as,
$$
f_{L_{\alpha,\theta}}(y)=Cy^{\alpha\theta-1}\int_{0}^{\infty}s^{-\theta
\alpha}{\mbox
e}^{-{(y/s)}^{\alpha}}f_{\alpha}(s)ds=Cy^{\alpha\theta-1}\int_{0}^{\infty}s^{-\theta
\alpha}\E[{\mbox e}^{-{(y/s)S_{\alpha}}}]f_{\alpha}(s)ds.
$$
for some constant C. Now it remains to write
$$
\E[{\mbox e}^{-{(y/s)S_{\alpha}}}]=\int_{0}^{\infty}{\mbox
e}^{-{vy/s}}f_{\alpha}(v)dv.
$$
The result is then obtained by algebraic manipulations. The second
equivalence is immediate since $L_{\alpha,\theta}$ is GGC$(\alpha
\theta,\varrho_{\alpha}).$ The remaining statements are
straightforward applications of beta-gamma calculus and
independence.\EndProof \Remark Statement [(iv)] generalizes the
known case of $G^{1/\alpha}_{1}\overset{d}=G_{1}/S_{\alpha}.$ It
is interesting to note that this special case can be interpreted
in terms of stochastic processes, as,
$$
\ell^{(\alpha)}_{G_{1}}\overset{d}=G^{\alpha}_{1}\frac{1}{{(S_{\alpha})}^{\alpha}}\overset{d}=G_{1}
$$
where $\ell^{(\alpha)}_{G_{1}}$ is the local time of a Bessel
process with dimension $2-2\alpha$ considered at the independent
exponential time $G_{1}.$  This description, and further
references, may be found in Chaumont and Yor~(2003, p. 114).  It
is also known from Barlow, Pitman and Yor~(1989) that the local
time for the bessel bridges,$\ell^{(\alpha,\alpha)}_{t}$ evaluated
at $G_{\alpha}$, satisfies
$\ell^{(\alpha,\alpha)}_{G_{\alpha}}\overset{d}=\ell^{(\alpha)}_{G_{1}}.$
Now using the scaling property in~\mref{scale} it follows from
statement[(iv)] that, $$ \ell^{(\alpha,\alpha\theta)}_{(G_{\theta
\alpha})}\overset{d}=\frac{{(G_{\theta
\alpha})}^{\alpha}}{{(T_{\alpha,\alpha\theta})}^{\alpha}}\overset{d}=G_{\theta}.
$$ See section 5 for an extension of this idea to more general
random times.\EndRemark

Theorem 4.1 establishes the distributional results,
$$
X_{\alpha}\overset{d}=B_{\alpha,1-\alpha}M_{\alpha}(\varrho_{\alpha})\overset{d}=M_{1}(\varrho^{(\alpha)}_{\alpha})=
{[(\ell^{(\alpha,\alpha)}_{B_{\alpha,1-\alpha}})]}^{1/\alpha}S_{\alpha}\overset{d}=
{[(\ell^{(\alpha)}_{S_{\alpha}})]}^{1/\alpha}
$$
where $(\ell^{(\alpha)}_{S_{\alpha}})$ denotes the local time of a
Bessel process evaluated at an independent unilateral
$(\alpha)$-stable random variable, $S_{\alpha}$. The fact that it
is representable as a Dirichlet mean functional,
$M_{1}(\varrho^{(\alpha)}_{\alpha})$ coupled with the explicit
density of $X_{\alpha},$ leads to an important identity below,
which plays a fundamental role in obtaining explicit densities
throughout the rest of this work.
\begin{prop}Define
$\Scr_{\alpha}(x)=\int_{0}^{\infty}\log(|x-y|)\varrho_{\alpha}(dy)=\E[\log(|x-X_{\alpha}|)].$
Then for $0<\alpha< 1,$
$$
\Scr_{\alpha}(x)=\frac{1}{2\alpha}\log(x^{2\alpha}+2x^{\alpha}\cos(\alpha
\pi)+1).
$$
Note that $\Scr_{\alpha}(x)$ is a special case of $\Phi(x).$
\end{prop}
\Proof From Theorem 3.1. it follows that the density of
$X_{\alpha}\overset{d}=M_{1}(\varrho^{(\alpha)}_{\alpha})$
satisfies the equivalence,
$$
\varrho_{\alpha}(x)=\frac{1}{\pi}\sin(\pi \alpha
[1-F_{X_{\alpha}}(x)]){\mbox e}^{-\Scr_{\alpha}(x)}x^{\alpha-1}.
$$
Solving this expressions for $\Scr_{\alpha}(x),$ and applying the
identity in~\mref{sinXid} concludes the result. \EndProof

 With this we obtain explicit expressions for the cdf
and density of $M_{\alpha \theta}(\varrho_{\alpha})$ as follows;
\begin{thm} From Theorem 4.1,
$M_{\alpha\theta}(\varrho_{\alpha})\overset{d}={[\ell^{(\alpha,\alpha\theta)}_{S_{\alpha}}]}^{1/\alpha}\overset{d}=S_{\alpha}/T_{\alpha,\alpha
\theta}.$ The form of the cdf for
$M_{\alpha\theta}(\varrho_{\alpha})$ for all $\alpha\theta>0$, is
given by~\mref{DPcdf}, with $\theta:=\alpha\theta,$ and
$$
\Delta_{
\alpha\theta}(x|\varrho_{\alpha})=\frac{1}{\pi}\frac{\sin(\pi
\theta \alpha
F_{X_{\alpha}}(x))}{{[x^{2\alpha}+2x^{\alpha}\cos(\alpha
\pi)+1]}^{{\theta}/{2}}}
$$
where $F_{X_{\alpha}}$ is given in~\mref{cdfXY}. Furthermore, a
general expression for the density is obtained
from~\mref{generaldensity} with $\theta:=\alpha\theta$ and \Eq
\label{genV}
d_{\alpha\theta}(x|\varrho_{\alpha})=\frac{\alpha\theta
x^{\alpha-1}}{\pi}\frac{[\sin(\pi \alpha[1-\theta
F_{X_{\alpha}}(x)])-x^{\alpha}\sin(\pi\theta \alpha
F_{X_{\alpha}}(x))]}{{[x^{2\alpha}+2x^{\alpha}\cos(\alpha
\pi)+1]}^{{\theta}/{2}+1}}.\EndEq\qed
\end{thm}
\Proof The bulk of the result is a straightforward application of
Proposition 4.2 combined with the explicit forms of the densities
in Section 3. The expression in~\mref{genV} follows by, noting the
explicit form of the density $\varrho_{\alpha}$, and applying the
identity $\sin(w-z)=\sin(w)cos(z)-\sin(z)cos(w),$ with $w=\pi
\alpha$ and $z=\pi\theta \alpha F_{X_{\alpha}}(x)$ \EndProof

From this we obtain the most explicit case of $\alpha\theta=1$,
and the important case where $\alpha\theta=\alpha$, related to the
local time of a Bessel Bridge evaluated at an independent stable
time.
\begin{cor}Consider the random variables in Theorem 4.2 then,
for $\theta=1$, corresponding to the Bessel bridge, the random
variable
$M_{\alpha}(\varrho_{\alpha})\overset{d}={[\ell^{(\alpha,\alpha)}_{S_{\alpha}}]}^{1/\alpha}\overset{d}=S_{\alpha}/T_{\alpha,\alpha}$
has density determined by~\mref{generaldensity} where
$d_{\alpha}(x|\varrho_{\alpha})$ simplifies to,
$$
d_{\alpha}(x|\varrho_{\alpha})=\frac{\alpha
x^{\alpha-1}}{\pi}\frac{(1-x^{2\alpha})\sin(\pi
\alpha)}{{[x^{2\alpha}+2x^{\alpha}\cos(\alpha \pi)+1]}^{2}}{\mbox
{ for }} x>0.
$$
When $\theta={1/\alpha}$, the density of
$M_{1}(\varrho_{\alpha})\overset{d}={[\ell^{(\alpha,1)}_{S_{\alpha}}]}^{1/\alpha}\overset{d}=S_{\alpha}/T_{\alpha,1}$
is given by
$$
f_{M_{1}}(x|{\varrho_{\alpha}})=\Delta_{1}(x|\varrho_{\alpha})=\frac{1}{\pi}\frac{\sin(\pi
F_{X_{\alpha}}(x))}{{[x^{2\alpha}+2x^{\alpha}\cos(\alpha
\pi)+1]}^{1/{2\alpha}}}.
$$
\end{cor}
\Proof Apply ~\mref{sinXid} to obtain the expression for
$d_{\alpha}(x|\varrho_{\alpha}).$\EndProof \Remark It is evident
from Theorem 3.1 and Proposition 4.2,  that we can also obtain
similar types of expressions, as in Theorem 4.2, for the densities
of
$M_{\theta}(\varrho^{(\alpha)}_{\alpha})\overset{d}=B_{\alpha\theta,\theta(1-\alpha)}M_{\alpha
\theta }(\varrho_{\alpha}),$ for all $\theta>0.$ For brevity we do
not give that here as it can be easily deduced. Instead we
concentrate on the interesting class of $M_{1}(\varrho^{(\alpha
\theta)}_{\alpha})\overset{d}=B_{\alpha\theta,1-\alpha\theta}M_{\alpha\theta}(\varrho_{\alpha})$
for $0<\alpha \theta<1.$ These have the densities of non-integral
form. \EndRemark

\begin{prop}For $0<\alpha\theta<1$,
$$M_{1}(\varrho^{(\alpha\theta)}_{\alpha})\overset{d}=B_{\alpha\theta,1-\alpha\theta}M_{\alpha\theta}(\varrho_{\alpha})
\overset{d}={[(\ell^{(\alpha,\alpha\theta)}_{B_{\alpha\theta,1-
\alpha\theta}})]}^{1/\alpha}S_{\alpha}$$ and has density \Eq
\label{mden} \frac{1}{\pi}\frac{\sin( \theta
C_{\alpha}(x))x^{\alpha\theta-1}}{{[x^{2\alpha}+2x^{\alpha}\cos(\alpha
\pi)+1]}^{{\theta}/{2}}} \EndEq where
$C_{\alpha}(x)=\cot^{-1}\left(\cot(\pi
\alpha)+\frac{x^{\alpha}}{\sin(\pi \alpha)}\right).$ The density
equates with $\varrho_{\alpha}$ when $\theta=1.$ For $0<\theta<1$
~\mref{mden} equates with the density of $B^{1/\alpha}_{\theta,1-
\theta}X_{\alpha}.$
\end{prop}
\Proof This is a consequence of Theorems 3.1 and 4.1 \EndProof

When $\alpha=1/k$, we obtain representations in terms of products
of gamma random variables.

\begin{prop}
Suppose that $\alpha=1/k$ for some integer $k=2,3\ldots$. Then
since for $\theta>0,$
$L_{1/k,\theta}\overset{d}=S_{1/k}G^{k}_{\theta}\overset{d}=G_{\theta/k}M_{\theta/k}(\varrho_{1/k})$,
one deduces that,
$$
M_{\theta/k}(\varrho_{1/k})=\frac{S_{1/k}}{T_{1/k,\theta/k}}\overset{d}=k^{k}S_{1/k}\prod_{j=1}^{k-1}G_{\frac{\theta+j}{k}}
$$
This implies
 $T_{1/k,\theta/k},$ is such that
$$
\frac{1}{T_{1/k,\theta/k}}\overset{d}=k^{k}\prod_{j=1}^{k-1}G_{\frac{\theta+j}{k}}
$$ distribution. One may obtain further representations by using
the known fact[see Chaumont and Yor~(2003, p. 113), that
$1/S_{1/k}\overset{d}=k^{k}\prod_{j=1}^{k-1}G_{j/k}.$
\end{prop}
\Proof The result may be deduced from the equivalence
$G^{k}_{\theta}S_{1/k}=G_{\theta/k}M_{\theta/k}(\varrho_{1/k})$
established in Theorem 4.1 and the following identity
$$
G^{k}_{\theta}={(G_{k(\frac{\theta}{k})})}^{k}=k^{k}G_{\frac{\theta}{k}}\prod_{j=1}^{k-1}G_{(\theta+j)/k},
$$
which is found in Chaumont and Yor~(2003, p. 113). \EndProof
\Remark As mentioned in Section 3, it is known from Cifarelli and
Melilli~(2000, p. 1394) that
$M_{\theta/2}(\varrho_{1/2})\overset{d}=G_{(\theta+1)/2}/G_{1/2},$
for $\theta\ge 2$. However the other results in Proposition 4.4
are apparently new for integer values, $k> 2$, and all $\theta>0$.
\EndRemark
\subsection{Results for tilted Linnik laws and occupation times}
We now address the case of the tilted processes which connects
with the occupation times.
\begin{prop} Suppose that $Z$ is GGC$^{(b,c)}(\alpha\theta,\varrho_{\alpha})$ Then, its L\'evy measure may be read from Propositions 3.1. and 3.2
In particular a GGC$^{(c,c)}(\alpha\theta,\varrho)$ is a
GGC$(\alpha\theta,\Lambda_{\alpha,p})$ with L\'evy measure
expressible as
$$
\alpha \theta s^{-1}\E[{\mbox e}^{-s/A_{\alpha,p}}]
$$
Additionally for $L_{c}$ a GGC$(\alpha\theta,\Lambda_{\alpha,p})$
random variable, we have the following distributional equivalences
\begin{itemize}\item[(i)] For all $\theta>0,$
$ L_{c}\overset{d}=G_{\theta
\alpha}P_{\alpha,\alpha\theta}(C)\overset{d}=G_{\theta
\alpha}M_{\alpha \theta}(\Lambda_{\alpha,p})$, hence $M_{\alpha
\theta}(\Lambda_{\alpha,p})\overset{d}=P_{\alpha,\theta}(C)$
\item[(ii)]As a special case, when $\theta=1$,
$A^{br}_{\alpha,p}\overset{d}=M_{\alpha}(\Lambda_{\alpha,p})$
\item[(iii)]When $\theta=1$,
$
L_{c}\overset{d}=G_{1}A^{*}_{\alpha,p}\overset{d}=G_{1}B_{\alpha,1-\alpha}M_{\alpha}(\Lambda_{\alpha,p})\overset{d}=G_{1}B_{\alpha,1-\alpha}P_{\alpha,\alpha}(C)$,
where
$A^{*}_{\alpha,p}\overset{d}=B_{\alpha,1-\alpha}M_{\alpha}(\Lambda_{\alpha,p})\overset{d}=B_{\alpha,1-\alpha}P_{\alpha,\alpha}(C)\overset{d}=G_{1}M_{1}(\Lambda^{(\alpha)}_{\alpha,p}),$
has density \Eq \label{gden}
\frac{(1-x)}{1-p}\Lambda_{\alpha,p}(dx) \EndEq Note
$A^{*}_{\alpha,p}$ equates with the random variable of Pitman and
Yor(1997b, Proposition 15).
\end{itemize}
\end{prop}
\Proof The correspondence to the $P_{\alpha,\alpha\theta}(C)$
models follows from~\mref{CSBesBridge} as mentioned earlier. All
other results are consequences of Proposition 3.2. and the use of
Beta Gamma calculus \EndProof

\Remark From Pitman and Yor~(1997b) and Barlow, Pitman and
Yor~(1989, p. 307),  let $g_{\alpha}=\sup\{t:t\leq 1:
B^{(\alpha)}_{t}=0\}$, where it is known that
$g_{\alpha}\overset{d}=B_{\alpha,1-\alpha}.$ Proposition 15 in
Pitman and Yor~(1997b) and their subsequent discussion establish
the results
$A_{\alpha,p}(g_{\alpha})\overset{d}=g_{\alpha}A^{br}_{\alpha,p}.$
with density ~\mref{gden}. Our result shows this random variable
is equivalent in distribution to the  Dirichlet mean functional,
$M_{1}(\Lambda^{(\alpha)}_{\alpha,p})\overset{d}=B_{\alpha,1-\alpha}M_{\alpha}(\Lambda_{\alpha,p}).$
Note that in view of Proposition 3.2. this can be seen to arise by
first noting again that
$M_{1}(\varrho^{(\alpha)}_{\alpha})\overset{d}=B_{\alpha,1-\alpha}M_{\alpha}(\varrho_{\alpha})\overset{d}=X_{\alpha}$
where $cX_{\alpha}/(cX_{\alpha}+1)\overset{d}=A_{\alpha,p}$ and
then applying the change of measure to get the distribution of
$M_{1}(\Lambda^{(\alpha)}_{\alpha,p})$. As mentioned in section 3,
the construction of Theorem 3.1, and indeed Theorem 3.2, were
devised in part to understand the mechanics of the result of
Pitman and Yor~(1997b) in more generality. The next result extends
this idea.\EndRemark

As mentioned in the previous remark the next result provides an
extension of the distributional result of Pitman and Yor~(1997b,
Proposition 15).  It also will yield our first concrete example of
the \emph{occupation time} of a  randomly skewed Bessel process.
This may also be seen as a precursor to the quantities we will
encounter in Section 5.

\begin{prop}For $0<\theta<1,$ and $G_{\theta}$ and $G_{1-\theta}$
independent, define
$\xi^{*}_{\theta}\overset{d}={pB_{\theta,1-\theta}}/({[q+pB_{\theta,1-\theta}]})\overset{d}={pG_{\theta}}/({G_{\theta}+G_{1-\theta}q}).$
Then define \Eq \label{FirstA}
A_{\alpha,\xi^{*}_{\theta}}\overset{d}=\frac{cB^{1/\alpha}_{\theta,1-\theta}X_{\alpha}}{cB^{1/\alpha}_{\theta,1-\theta}X_{\alpha}+1}=\frac{cM_{1}(\varrho^{(\alpha\theta)}_{\alpha})}{cM_{1}(\varrho^{(\alpha\theta)}_{\alpha})+1}
\EndEq  The explicit density of the quantity in~\mref{FirstA} is
given by,
$$
f_{\alpha,\theta}(y)=
\frac{q^{\theta}}{\pi}\frac{y^{\alpha\theta-1}{(1-y)}^{-1}\sin(\theta
C_{\alpha}(\frac{q^{1/\alpha}y}{p^{1/\alpha}(1-y)}))}{{[y^{2\alpha}q^{2}+2qpy^{\alpha}{(1-y)}^{\alpha}\cos(\alpha
\pi)+{(1-y)}^{2\alpha}p^{2}]}^{{\theta}/{2}}}.
$$
The  density of $B_{\alpha\theta,1-\alpha \theta}P_{\alpha,\alpha
\theta}(C)\overset{d}=B_{\alpha\theta,1-\alpha
\theta}A^{(\alpha,\alpha \theta)}_{\alpha,p}\overset{d}=
B_{\alpha\theta,1-\alpha \theta} M_{\alpha
\theta}(\Lambda_{\alpha,p})\overset{d}=M_{1}(\Lambda^{(\alpha\theta)}_{\alpha,p})$,
is \Eq \frac{1-y}{{(1-p)}^{\theta}}f_{\alpha,\theta}(y)=
\frac{1}{\pi}\frac{y^{\alpha\theta-1}\sin(\theta
C_{\alpha}(\frac{q^{1/\alpha}y}{p^{1/\alpha}(1-y)}))}{{[y^{2\alpha}q^{2}+2qpy^{\alpha}{(1-y)}^{\alpha}\cos(\alpha
\pi)+{(1-y)}^{2\alpha}p^{2}]}^{{\theta}/{2}}} \label{Bden}.\EndEq
More generally the results hold for all $0<\theta \alpha<1,$ where
$f_{\alpha,\theta}(y)$ is the density of the random variable
$cM_{1}(\varrho^{(\alpha\theta)}_{\alpha})/(cM_{1}(\varrho^{(\alpha\theta)}_{\alpha})+1),$
and~\mref{Bden} is the density of
$M_{1}(\Lambda^{(\alpha\theta)}_{\alpha,p}).$
\end{prop}
\Proof The density of~\mref{FirstA} is obtained from Proposition
4.3. To obtain,~\mref{Bden} note that we are using the fact
$G_{\alpha\theta}P_{\alpha
,\alpha\theta}(C)\overset{d}=G_{1}B_{\alpha\theta,1-\alpha\theta}P_{\alpha,\alpha\theta}(C)$.
Now it becomes evident that this distribution is obtained from
exponentially tilting
${(p/q)}^{1/\alpha}G_{1}M_{1}(\varrho^{(\alpha\theta)}_{\alpha})$
in the sense of Proposition 3.2 and Proposition 4.3.
 \EndProof
\Remark At this point we could use Theorem 4.2, combined with
Proposition 3.2, to obtain expressions for the density and cdf of
$A^{(\alpha,\alpha\theta)}_{\alpha,p}\overset{d}=P_{\alpha,\alpha\theta}(C)=M_{\alpha\theta}(\Lambda_{\alpha,p})$.
Or one could use directly ~\mref{invA} and~\mref{cdfA}. These
would provide alternative expressions for
$P_{\alpha,\alpha\theta}(C)$ obtained in James, Lijoi and
Pr\"unster~(2006). In that work, the authors addressed the case of
more general functionals $P_{\alpha,\alpha\theta}(g)$, where they
obtained these laws by a direct inversion of the appropriate
Cauchy-Stieltjes transform. That work does not address functionals
such as $M_{\alpha\theta}(\varrho_{\alpha})$.  A description of
the laws of $A^{(\alpha,\alpha\theta)}_{\alpha,p}$ will appear as
a special case of the forthcoming Proposition 5.11.\EndRemark

\section{GGC/FGGC stable compositions and occupation laws for randomly skewed Bessel processes} Raising things to the level
of processes we recall the known fact that a positive Linnik
process equates to
$L_{\alpha,\theta}(t)\overset{d}=S_{\alpha}(\Gamma_{\theta}(t))\overset{d}=G^{1/\alpha}_{\theta
t}S_{\alpha}$, where in the second equality $S_{\alpha}(t)$ is a
stable subordinator independent of the gamma subordinator,
$\Gamma_{\theta}(t).$ Based on our previous results we now study
the class of models
$S_{\alpha}(Z_{\theta}(t))\overset{d}={[Z_{\theta}(t)]}^{1/\alpha}S_{\alpha}$
where $Z_{\theta}$ is a GGC/FGGC. A special case is where
$Z_{\theta}(t)$ is itself a stable subordinator of index
$0<\beta<1,$ say $S_{\beta}(t)$, which satisfies the  important
identity
$S_{\alpha}(S_{\beta}(1))\overset{d}={[S_{\beta}]}^{1/\alpha}S_{\alpha}\overset{d}=S_{\alpha\beta},$
an $(\alpha\beta)$-stable random variable. Note that in this case
$Z_{\theta}=S_{\beta}$ is not a FGGC. \Remark Obviously
compositions of a stable subordinator have been previously studied
from several important perspectives. For instance, this operation
has recently been shown to play an interesting role in
applications involving coagulation/fragmentation phenomena as
described in Pitman~(2006) and Bertoin~(2006).\EndRemark

Here, using the property that $Z_{\theta}(t)\overset{d}=G_{\theta
t}M_{\theta t}(\nu),$ we take a different view of
$S_{\alpha}(Z_{\theta}(t))$ as being equivalent in distribution to
$G^{1/\alpha}_{\theta t}S_{\alpha}(t){[M_{\theta
t}(\nu)]}^{1/\alpha}$. That is to say, the viewpoint that these
are scale mixtures of positive Linnik random variables, which
leads to a variety of interesting consequences. Lets call this
class GGC$_{\alpha}(\theta,\nu)$ which will also denote the law of
the random variable $S_{\alpha}(Z_{\theta}(1)).$

Now, from the description of the L\'evy density of a Linnik random
variable in Section 4, it is evident that the Laplace transform,
at $t=1$, of $S_{\alpha}(Z_{\theta}(1)),$ can be expressed as
$$
\E[{(1+\lambda^{\alpha}M_{\theta}(\nu))}^{-\theta}]={\mbox
e}^{-\psi^{(\alpha)}_{\theta}(\lambda)}
$$
where
$$
\psi^{(\alpha)}_{\theta}(\lambda)=\int_{0}^{\infty}\theta\log(1+\lambda^{\alpha}r)\nu(dr)=
\alpha\theta\int_{0}^{\infty}\E[\log(1+\lambda
X_{\alpha}r^{1/\alpha})]\nu(dr).
$$

That is, the L\'evy measure of
$S_{\alpha}(Z_{\theta}(1))\overset{d}=G^{1/\alpha}_{\theta}S_{\alpha}{[M_{\theta}(\nu)]}^{1/\alpha}$
can be expressed as
$$
\alpha\theta s^{-1}\int_{0}^{\infty}\E[{\mbox
e}^{-s/(X_{\alpha}r^{1/\alpha})}]\nu(dr){\mbox { for }}s>0.
$$

\subsection{Tilting and Randomly Skewed Bessel Bridges} In parallel to Section 4,  we also discuss its
tilted version, which as we shall see connects naturally with the
idea of randomly skewed Bessel processes and bridges. Using the
same type of exponential tilting as in Section 4 we call the
resulting random variables $GGC^{(b,c)}_{\alpha}(\theta,\nu).$
Again, provided that $c\neq 0,$ we see that if a random variable
$L$ is $GGC^{(b,c)}_{\alpha}(\theta,\nu)$, then
$L\overset{d}=(b/c)L^{*}$ where $L^{*}$ is
$GGC^{(c,c)}_{\alpha}(\theta,\nu).$ An initial description of the
Laplace transform of  $L^{*}$ is \Eq \label{stabletilt}
\frac{\E[{(1+{(1+\lambda)}^{\alpha}c^{\alpha}M_{\theta}(\nu))}^{-\theta}]}{\E[{(1+c^{\alpha}M_{\theta}(\nu))}^{-\theta}]}
={\mbox e}^{-\Psi^{(\alpha)}_{\theta}(\lambda)} \EndEq where $
\Psi^{(\alpha)}_{\theta}(\lambda)=\psi^{(\alpha)}_{\theta}(c(1+\lambda))-\psi^{(\alpha)}_{\theta}(c)$.
Now, by setting $p_{\alpha}(r)=c^{\alpha}r/1+c^{\alpha}r)$ and
algebra, this is equivalent to \Eq \label{LevyAA}
\Psi^{(\alpha)}_{\theta}(\lambda)=\alpha\theta\int_{0}^{\infty}\E[\log(1+\lambda
A_{\alpha,p_{\alpha}(r)})]\nu(dr)=\alpha\theta\int_{0}^{\infty}\log(1+\lambda
x)\Lambda_{\alpha,\nu}(dx) \EndEq where
$$
\Lambda_{\alpha,\nu}(dx)=\int_{0}^{\infty}\Lambda_{\alpha,p_{\alpha}(r)}(dx)\nu(dr).
$$
 Before examining this quantity further we first formally define,
albeit briefly,  what we mean by randomly skewed Bessel processes
and bridges and their occupation times. Let $\xi$ denote a random
variable on $[0,1]$ chosen independently of a Bessel process
$\{B^{(\alpha)}(t):t\ge 0\}$. Then the process
$B^{(\alpha)}_{\xi}(t)$ is said to be a $\xi$-randomly skewed
Bessel process if $B^{(\alpha)}_{\xi}|\xi=p$ is a $p$-skewed
Bessel process as defined in Barlow, Pitman and Yor~(1989). That
is, $\Pe(B^{(\alpha)}_{\xi}(s)>0|\xi)=\xi.$  Say that
$A_{\alpha,\xi}(t):=\int_{0}^{t}\indic(B^{(\alpha)}_{\xi}(s)>0)ds$
is the time spent positive of a $\xi$-randomly skewed Bessel
process up to time $t.$ Then conditional on $\xi=p$, it has
distribution $A_{\alpha,p}(t)$ . Equivalently, one has, for
$A_{\alpha,\xi}:=A_{\alpha,\xi}(1)$, \Eq \label{Arep}
A_{\alpha,\xi}\overset{d}=\frac{\xi^{1/\alpha}X_{\alpha}}{\xi^{1/\alpha}X_{\alpha}+{(1-\xi)}^{1/\alpha}}.
\EndEq Based on \mref{Arep} we present an interesting special
case.
\begin{prop}Define $p_{\alpha}:=p^{\alpha}/[p^{\alpha}+q^{\alpha}].$ For $0<\beta<1$, let
$\xi\overset{d}=A_{\beta,p_{\alpha}}=c^{\alpha}X_{\beta}/(c^{\alpha}X_{\beta}+1)$
which is equivalent to $\xi/(1-\xi)=c^{\alpha}X_{\beta}$. Then
$$
A_{\alpha,\xi}\overset{d}=A_{\alpha\beta,p}.
$$
\end{prop}
\Proof The result follows from the definition in~\mref{Arep} and
the easily verified fact that
$[X_{\beta}]^{1/\alpha}X_{\alpha}\overset{d}=X_{\alpha\beta}$
\EndProof
 $\xi$-randomly skewed bridges are defined in an analogous
manner where conditionally on $\xi=p$ they are $p$-skewed Bessel
bridges. Denote the random variable corresponding to the time
spent positive up till time $1$ of such a process as
$A^{br}_{\alpha,\xi}$, which conditional on $\xi=p$, is equivalent
in distribution to
$A^{br}_{\alpha,p}\overset{d}=P_{\alpha,\alpha}(C),$ having the
Cauchy-Stieltjes transform of order $\alpha\theta$
in~\mref{CSBesBridge}. More generally, by the usual
PD$(\alpha,\alpha\theta)$ change of measure, we can define
$A^{(\alpha,\alpha\theta)}_{\alpha,\xi}$ which conditionally on
$\xi=p$ equates in distribution with
$A^{(\alpha,\alpha\theta)}_{\alpha,p}\overset{d}=P_{\alpha,\alpha\theta}(C).$
Now we shall describe a sub-class of such processes which equates
with a random variable described by the
transform~\mref{stabletilt}.
\begin{prop}Let
$A^{(\alpha,\alpha\theta)}_{\alpha,\xi_{\theta}}$ correspond to a
$\xi_{\theta}$-skewed occupation time with $\xi_{\theta}$ having
the specific density
$Q_{\theta}(du|\nu)/du=\kappa_{\theta}{(1-u)}^{\theta}\tilde{Q}_{\theta}(du|\nu)/du$
where
$\kappa_{\theta}=1/\E[{(1+c^{\alpha}M_{\theta}(\nu))}^{-\theta}]={\mbox
e}^{\psi^{(\alpha)}_{\theta}(c)}$ and
$$
\tilde{Q}_{\theta}(du)=\frac{1}{{(1-u)}^{2}}f_{c^{\alpha}M_{\theta}}(\frac{u}{1-u}|\nu)du
$$
That is $\tilde{Q}_{\theta}$ is the distribution of the random
variable
$c^{\alpha}M_{\theta}(\nu)/(c^{\alpha}M_{\theta}(\nu)+1).$
Equivalently
$\xi_{\theta}\overset{d}=M_{\theta}(\nu^{(c^{\alpha},c^{\alpha})})$
and  $A^{(\alpha,\alpha\theta)}_{\alpha,\xi_{\theta}}$ satisfies
the following properties.
\begin{enumerate}
\item[(i)]The Cauchy-Stieltjes transform of order $\alpha \theta$ of
$A^{(\alpha,\alpha\theta)}_{\alpha,\xi_{\theta}}$ is expressible
as
$$
\int_{0}^{1}\E[{(1+\lambda
A^{(\alpha,\alpha\theta)}_{\alpha,p})}^{-\alpha\theta}]Q_{\theta}(dp|\nu)=\int_{0}^{1}{(q+p{(1+\lambda)}^{\alpha})}^{-\theta}Q_{\theta}(dp|\nu)
$$
and equals~\mref{stabletilt}.
\item[(ii)]Hence, by ~\mref{LevyAA},
$A^{(\alpha,\alpha\theta)}_{\alpha,\xi_{\theta}}\overset{d}=M_{\alpha\theta}(\Lambda_{\alpha,\nu}),$
and
$L\overset{d}=G_{\alpha\theta}A^{(\alpha,\alpha\theta)}_{\alpha,\xi_{\theta}}$
is GGC$(\alpha\theta,\Lambda_{\alpha,\nu}).$
\end{enumerate}
\end{prop}
\Proof This result is easily verified by an argument similar to
the proof of Theorem 3.2. Hence we omit the details.\EndProof

Now define \Eq \label{genvar}
\varrho_{\alpha,\nu}(x)=\frac{\sin(\pi\alpha)}{\pi}\int_{0}^{\infty}
\frac{
x^{\alpha-1}r}{x^{2\alpha}+2x^{\alpha}r\cos(\alpha\pi)+r^{2}}\nu(dr).
\EndEq
 We can summarize these results in an equivalent manner;
\begin{prop}Let $L\overset{d}={[Z_{\theta}(1)]}^{1/\alpha}S_{\alpha}$ denote a GGC$_{\alpha}(\theta,\nu)$ random variable
then \Eq\label{idtime}L\overset{d}=G_{\theta
\alpha}M_{\alpha\theta}(\varrho){[M_{\theta}(\nu)]}^{1/\alpha}=G_{\alpha\theta}M_{\alpha\theta}(\varrho_{\alpha,\nu}).\EndEq
where $\varrho_{\alpha,\nu}$ is defined in~\mref{genvar}. That is
GGC$_{\alpha}(\theta,\nu)$ is equivalent to
GGC$(\alpha\theta,\varrho_{\alpha,\nu}).$ Now without loss of
generality if $c\neq 0$, then set $c=1$ and consider the case
where ${\tilde L}$ is GGC$_{\alpha}(\theta,\nu^{(1,1)})$ with
Laplace transform as in~\mref{stabletilt}. Then ${\tilde
L}\overset{d}=G_{\alpha\theta}A^{(\alpha,\alpha\theta)}_{\alpha,\xi_{\theta}}\overset{d}=G_{\alpha\theta}M_{\alpha\theta}(\Lambda_{\alpha,\nu})$
is GGC$(\alpha \theta,\Lambda_{\alpha,\nu})$ In general if $L$ is
GGC$^{(b,c)}_{\alpha}(\theta,\nu)$ random variable, then its
L\'evy measure can be expressed as,
$$
\alpha \theta s^{-1}{\mbox e}^{-(c/b)s}\int_{0}^{\infty}\E[{\mbox
e}^{-s/[bX_{\alpha}r^{1/\alpha}]}]\nu(dr).\qed
$$
\end{prop}
The next result is specialized to the case where $\nu=H.$
\begin{prop} Suppose that $L$ is
GGC$(\alpha\theta,\varrho_{\alpha,H})$, then there exists a random
variable $R$ with distribution $H$, such that
$R^{1/\alpha}X_{\alpha}$ has density $\varrho_{\alpha,H}$ and L is
a FGGC with L\'evy measure
$$
\alpha\theta s^{-1}\E[{\mbox e}^{-s/[R^{1/\alpha}X_{\alpha}]}].
$$
Correspondingly ${\tilde L}$, the exponential tilt of $L$ with
$c=1$, is GGC$(\alpha\theta,\Lambda_{\alpha,H})$ and has the
L\'evy measure
$$
\alpha\theta s^{-1}\E[{\mbox e}^{-s/A_{\alpha,\xi}}]
$$
where $\xi/(1-\xi)\overset{d}=R$ and
$A_{\alpha,\xi}\overset{d}=R^{1/\alpha}X_{\alpha}/(R^{1/\alpha}X_{\alpha}+1)$
has distribution $\Lambda_{\alpha,H}.$ Hence in this setting
$A^{(\alpha,\alpha\theta)}_{\alpha,\xi_{\theta}}\overset{d}=M_{\alpha\theta}(\Lambda_{\alpha,H}),$
and $\xi_{\theta}=M_{\theta}(H^{(c^{\alpha},c^{\alpha})})$ are
Dirichlet mean functionals.\qed
\end{prop}

 \Remark
Bondesson~(1992) describes various features of the class which we
call GGC$_{\alpha}(\theta,\nu)$. See for instance Theorem 3.3.2 of
that work which establishes the fact that this class of models are
indeed GGC. As noted by Bondesson~(1992), compositions of GGC
random variables with some GGC subordinator are not always GGC. An
example is the composition of a Gamma subordinator with another
GGC process. Note that our representation of the L\'evy measure in
terms of $X_{\alpha}R^{1/\alpha}$ and our subsequent usage of it
appears to be new. \EndRemark \Remark It is important to note that
special cases of such models have already appeared in the
literature. Proposition 19 of Pitman~(1999), in connection with
coagulation phenomena, shows, with obvious rephrasing,  that the
law of $A^{br}_{\alpha,\xi}$, where $\xi=B_{\beta,1-\beta}$,
corresponds to a \textsc{Beta}$(\alpha-\alpha\beta, \alpha\beta)$
random variable. Another interesting case where a model having a
particular distribution of $A^{br}_{\alpha,\xi}$ may be found is
in Aldous and Pitman~(2004). Specifically, the distribution of a
random lengths $1-T_{k}$ of a \emph{T-partition }, as described in
equation (67) of Aldous and Pitman~(2004), may be re-expressed as
$$\Pe(T_{k}\in dt)=\int_{0}^{1}\Pe(A^{br}_{\alpha,p}\in dx)\Pi_{k}(dp)=\Pe(A^{br}_{\alpha,\xi}\in dx)$$
where $\Pi_{k}(dp)/dp={(-\log(p))}^{k-1}/(k-1)!$  and hence
$\xi\overset{d}=\prod_{i=1}^{k}U_{i},$ for $(U_{i})$ independent
Uniform$[0,1]$ random variables. We now show how to use a
construction of Aldous and Pitman~(2004) that equates, more
generally, with the distribution of
$A^{(\alpha,\alpha\theta)}_{\alpha,\xi}.$ Extending the definition
in Aldous and Pitman~(2004, p.24), define for $0<u<1$
$$
\tau^{(\alpha,\alpha\theta)}_{\alpha,u}=\inf\{t:\frac{\ell^{(\alpha,\alpha\theta)}_{t}}{\ell^{(\alpha,\alpha\theta)}_{1}}=u\}.
$$
Then it follows that
$\tau^{(\alpha,\alpha\theta)}_{\alpha,u}\overset{d}=A^{(\alpha,\alpha\theta)}_{\alpha,u}$
and hence, $
\tau^{(\alpha,\alpha\theta)}_{\alpha,\xi}\overset{d}=A^{(\alpha,\alpha\theta)}_{\alpha,\xi}.
$ \EndRemark

\Remark We note that one could use the results for the laws of
$A^{(\alpha,\alpha\theta)}_{\alpha,p}\overset{d}=P_{\alpha,\alpha\theta}(C)$
in James, Lijoi and Pr\"unster~(2006), or the form obtainable from
Section 4, to represent the laws of
$A^{(\alpha,\alpha\theta)}_{\alpha,\xi}$ by mixing with respect to
a distribution of $\xi.$ However it is not a trivial matter to
simplify such expressions. Here we will use a direct approach.
\EndRemark

\subsection{Random time changes for Local and Occupation
times} Here we translate the results in the previous section in
terms of random time changes of the local and occupation times.
\begin{prop} Let $(\ell^{(\alpha,\alpha\theta)}_{t}; t\ge 0)$ denote the local time under
the PD$(\alpha,\alpha\theta )$law for $\alpha\theta>0$. Let
$Z_{\theta}(1)$ denote a GGC$(\theta,\nu)$ random variable. Then
it follows that
\begin{enumerate}
\item[(i)]
$\ell^{(\alpha,\alpha\theta)}_{G_{\theta
\alpha}M^{1/\alpha}_{\theta}(\nu)}\overset{d}=\ell^{(\alpha,\alpha\theta)}_{G_{\theta
\alpha}}M_{\theta}(\nu)\overset{d}=Z_{\theta}(1)\overset{d}=G_{\theta}M_{\theta}(\nu)$
\item[(ii)]
$
A_{\alpha,p}(\tau_{\alpha}(Z_{\theta}(1)))\overset{d}=p^{1/\alpha}{[(\ell^{(\alpha,\theta)}_{G_{\theta
\alpha}})]}^{1/\alpha}S_{\alpha}M^{1/\alpha}_{\theta}(\nu)\overset{d}=p^{1/\alpha}\tilde{L}_{\alpha,\theta}.
$
\end{enumerate}
\end{prop}
\Proof Statement [(i)] follows from the scaling property
in~\mref{scale}, ~\mref{idtime}, and the fact that
$M_{\alpha\theta}(\varrho_{\alpha})\overset{d}={[\ell^{(\alpha,\alpha\theta)}_{S_{\alpha}}]}^{1/\alpha}.$
Statement [(ii)] is straightforward. \EndProof
 The result states that we can arrange for the time changed
local times to have any GGC distribution. We present the next
result for a simple illustration
\begin{prop} Let $(\ell^{(\alpha,\alpha\theta)}_{t}; t\ge 0)$ denote the local time under
the PD$(\alpha,\alpha\theta )$law for $\alpha\theta>0$. Set
$M_{\theta}(\nu)=B_{\alpha\theta,\theta(1-\alpha)}.$ Then it
follows that
\begin{enumerate}
\item[(i)]
$\ell^{(\alpha,\alpha\theta)}_{G_{\theta
\alpha}B^{1/\alpha}_{\alpha\theta,\theta(1-\alpha)}}=\ell^{(\alpha,\alpha\theta)}_{G_{\theta
\alpha}}B_{\alpha\theta,\theta(1-\alpha)}\overset{d}=G_{\alpha\theta}$
\item[(ii)]
$
A_{\alpha,p}(\tau_{\alpha}(G_{\alpha\theta}))\overset{d}=p^{1/\alpha}G^{1/\alpha}_{\alpha\theta}S_{\alpha}
$
\end{enumerate}
\end{prop}

\subsection{Scaling Calculus for GGC$_{\alpha}(\theta,\nu)$ mean
functionals} We now show that our results establish a rather
remarkable type of scaling calculus for mean functionals.
\begin{prop}If $L$ is GGC$_{\alpha}^{(b,c)}(\theta,\nu)$, it follows that
GGC$_{\alpha}^{(1,0)}(\theta,\nu)=GGC(\alpha
\theta,\varrho_{\alpha,\nu})$ and
GGC$_{\alpha}^{(1,1)}(\theta,\nu)=GGC(\alpha
\theta,\Lambda_{\alpha,\nu}).$ That is in the first case
$L\overset{d}=G_{\alpha\theta}M_{\alpha
\theta}(\varrho_{\alpha,\nu})$ and in the second
$L\overset{d}=G_{\alpha\theta}M_{\alpha
\theta}(\Lambda_{\alpha,\nu}).$ These points lead to the following
results \begin{enumerate}
\item[(i)]$M_{\alpha
\theta}(\varrho_{\alpha,\nu})\overset{d}=M_{\alpha
\theta}(\varrho_{\alpha}){[M_{\theta}(\nu)]}^{1/\alpha}$
\item[(ii)]Following Theorem 3.2 and Proposition 3.1, $A^{(\alpha,\alpha\theta)}_{\alpha,
\xi_{\theta}}=M_{\alpha
\theta}(\Lambda_{\alpha,\nu})\overset{d}=Y_{\alpha \theta,1}$,
with $M\overset{d}=M_{\alpha \theta}(\varrho_{\alpha,\nu}).$
\end{enumerate}
The results apply for arbitrary Thorin measures $\nu$ but yield
results specific to Dirichlet mean functionals by setting
$\nu=H.$\qed
\end{prop}

A variation of this result is given below,

\begin{prop}Theorem 3.1 and Proposition 5.7 imply,
\begin{enumerate}
\item[(i)]$M_{1}(\varrho^{(\alpha\theta)}_{\alpha,H})=B_{\alpha\theta,1-\alpha\theta}M_{\alpha
\theta}(\varrho_{\alpha,H})\overset{d}=M_{1}(\varrho^{(\alpha\theta)}_{\alpha}){[M_{\theta}(H)]}^{1/\alpha}$
\item[(ii)]$M_{1}(\varrho^{(\alpha)}_{\alpha,H})=X_{\alpha}{[M_{1}(H)]}^{1/\alpha}$
\item[(ii)]$M_{1}(\Lambda^{(\alpha\theta)}_{\alpha,H})\overset{d}=B_{\alpha\theta,1-\alpha\theta}A^{(\alpha,\alpha\theta)}_{\alpha,\xi}.$
\end{enumerate}
Notice that when $0<\theta<1$,
$M_{1}(\varrho^{(\alpha\theta)}_{\alpha,H})\overset{d}=X_{\alpha}B^{1/\alpha}_{\theta,1-\theta}{[M_{\theta}(H)]}^{1/\alpha}$
but
$B^{1/\alpha}_{\theta,1-\theta}{[M_{\theta}(H)]}^{1/\alpha}={[M_{1}(H^{(\theta)})]}^{1/\alpha}$.
Hence this is a special case of statement [(ii)].\qed
\end{prop}

Notice now that
$$
\varrho_{\alpha,H}(x)=\frac{\sin(\pi\alpha)}{\pi}\int_{0}^{\infty}
\frac{
x^{\alpha-1}r}{x^{2\alpha}+2x^{\alpha}r\cos(\alpha\pi)+r^{2}}H(dr).
$$
The cdf of $X_{\alpha}R^{1/\alpha}$,
$\Upsilon_{\alpha}(x)=\Pe(X_{\alpha}R^{1/\alpha}\leq x),$ can be
represented as \Eq \label{cdfX}
\Upsilon_{\alpha}(x)=1-\frac{1}{\pi\alpha}\int_{0}^{\infty}\cot^{-1}\left(\cot(\pi
\alpha)+\frac{x^{\alpha}/r}{\sin(\pi \alpha)}\right)H(dr) \EndEq
which simplifies in the case of $\alpha=1/2$ to,
$$
\Upsilon_{1/2}(x)=\frac{2}{\pi}\int_{0}^{\infty}\arctan\left(\frac{\sqrt{x}}{r}\right)H(dr).
$$
\begin{prop}Let
$X_{\alpha}=S_{\alpha}/S'_{\alpha},$ having density~\mref{denX}.
Let $R$ be a random variable with distribution $H$. Then the
random variable $X_{\alpha}R^{1/\alpha}$ has density
$\varrho_{\alpha,H}$. Define
$\Scr^{H}_{\alpha}(x)=\E[\log(|x-X_{\alpha}R^{1/\alpha}|)\indic(x\neq
XR^{1/\alpha})].$ Then,
$$
\Scr^{H}_{\alpha}(x)=\frac{1}{2\alpha}\int_{0}^{\infty}\log(x^{2\alpha}+2x^{\alpha}r\cos(\alpha\pi)+r^{2})H(dr)
$$ with derivative $$
\sigma^{H}_{\alpha}(x)=\int_{0}^{\infty}\frac{ x^{2\alpha-1}+
x^{\alpha-1}r\cos(\alpha\pi)}{x^{2\alpha}+2x^{\alpha}r\cos(\alpha\pi)+r^{2}}H(dr).
$$
$
\Scr^{H}_{1/2}(x)=\int_{0}^{\infty}\log(x+r^{2})H(dr)=\int_{0}^{\infty}\log(1+r^{2}/x)H(dr)+\log(x).
$ \qed
\end{prop}
The next result shows that $\Scr^{H}_{\alpha}(x)$ may be seen as
the L\'evy exponent of certain FGGC models. This provides another
way for further simplification.
\begin{prop} Let $R$ be a random variable with distribution $H.$ Define
$W_{\alpha,x}\overset{d}=2Rx^{-\alpha}\cos(\pi
\alpha)+R^{2}x^{-2\alpha}$ and denote its distribution as
$H^{x}_{\alpha}.$ Then,
$$
{\mbox e}^{-\theta \alpha \Scr^{H}_{\alpha}(x)}=\E[{\mbox
e}^{-G_{\theta/2}M_{\theta/2}(H^{x}_{\alpha})}]x^{-\theta \alpha}
$$
When $\alpha=1/2$ let $H_{1/2}(y)=H(\sqrt{y})$ denote the cdf of
$R^{2}$ then
$$
{\mbox e}^{-\theta/2 \Scr^{H}_{1/2}(x)}=\E[{\mbox
e}^{-\frac{1}{x}G_{\theta/2}M_{\theta/2}(H_{1/2})}]x^{-\theta/2}.
$$\qed
\end{prop}

\subsection{Density formula}
With this we can obtain explicit expressions for the cdf and
density of $M_{\alpha \theta}(\varrho_{\alpha.H})$ as follows
\begin{thm}
Recall that
$M_{\alpha\theta}(\varrho_{\alpha,H})\overset{d}=M_{\alpha\theta}(\varrho_{\alpha}){[M_{\theta}(H)]}^{1/\alpha}.$
The form of the cdf for $M_{\alpha\theta}(\varrho_{\alpha,H})$ for
all $\alpha\theta>0$, is given by~\mref{DPcdf}, with
$\theta:=\alpha\theta,$ and
$$
\Delta_{\theta \alpha}(x|\varrho_{\alpha,H})=\frac{1}{\pi}\sin(\pi
\theta \alpha \Upsilon_{\alpha}(x)){\mbox e}^{-\theta
\alpha\Scr^{H}_{\alpha}(x)}
$$
where $\Upsilon_{\alpha}$ is given in~\mref{cdfX}.
$\Delta_{1}(x|\varrho_{\alpha,H})$ is the density of
$M_{1}(\varrho_{\alpha,H}).$ Furthermore, a general expression for
the density is obtained from~\mref{generaldensity} with
$\theta:=\alpha\theta$ and $$
d_{\alpha\theta}(x|\varrho_{\alpha,H})=\frac{\alpha\theta}{\pi}{\mbox
e}^{-\theta \alpha\Scr^{H}_{\alpha}(x)}\[\sin(\pi\alpha[1-\theta
\Upsilon_{\alpha}(x)])B^{H}_{\alpha,1}(x)-\sin(\pi \theta \alpha
\Upsilon_{\alpha}(x))B^{H}_{\alpha,2}(x)
\],$$ where$$B^{H}_{\alpha,1}(x)=\int_{0}^{\infty}\frac{
x^{\alpha-1}r}{x^{2\alpha}+2x^{\alpha}r\cos(\alpha\pi)+r^{2}}H(dr)$$
and
$$
B^{H}_{\alpha,2}(x)=\int_{0}^{\infty}\frac{
x^{2\alpha-1}}{x^{2\alpha}+2x^{\alpha}r\cos(\alpha\pi)+r^{2}}H(dr).
$$
\qed
\end{thm}
Now applying Theorem 3.1 we obtain,
\begin{thm}For $0<\alpha \theta <1$, the density of
$M_{1}(\varrho^{(\alpha\theta)}_{\alpha,H})=B_{\alpha\theta,1-\alpha\theta}M_{\alpha\theta}(\varrho_{\alpha}){[M_{\theta}(H)]}^{1/\alpha}$
is given by
$$
\frac{1}{\pi}\sin(\pi \theta \alpha
[1-\Upsilon_{\alpha}(x)]){\mbox e}^{-\theta
\alpha\Scr^{H}_{\alpha}(x)}x^{\alpha\theta-1}
$$
When $\theta=1,$ the density of
$M_{1}(\varrho^{(\alpha)}_{\alpha,H})\overset{d}=X_{\alpha}{[M_{1}(H)]}^{1/\alpha}$
can also be expressed as,
$$
\varrho_{\alpha,F_{M_{1}(H)}}(x)=\frac{\sin(\pi\alpha)}{\pi}\int_{0}^{\infty}
\frac{
x^{\alpha-1}r}{x^{2\alpha}+2x^{\alpha}r\cos(\alpha\pi)+r^{2}}f_{M_{1}}(r|H)dr.
$$\qed
\end{thm}
Then next results are for the $\xi$-skewed occupation times.
\begin{prop}Assume without loss of generality that $c=1$. Then, for all $\alpha\theta>0$, the density of
$M_{\alpha
\theta}(\Lambda_{\alpha,H})\overset{d}=A^{(\alpha,\alpha\theta)}_{\alpha,\xi_{\theta}}$
can be expressed as
$$
\kappa_{\theta}{(1-y)}^{\alpha\theta-2}f_{M_{\alpha\theta}}(\frac{y}{1-y}|\varrho_{\alpha,H})
$$
\end{prop}
\Proof This is an application of Proposition 3.2. See also
statement [(ii)] of Proposition 5.7. \EndProof The next result is
a generalization of Proposition 4.6.
\begin{prop}Assume without loss of generality that $c=1$. Then, for all $0<\alpha\theta<1$, the density of
$M_{1}(\Lambda^{(\alpha\theta)}_{\alpha,H})\overset{d}=B_{\alpha\theta,1-\alpha\theta}A^{(\alpha,\alpha\theta)}_{\alpha,\xi_{\theta}}$
can be expressed as
$$
\kappa_{\theta} \frac{1}{\pi}\sin(\pi \theta \alpha
[1-\Upsilon_{\alpha}(\frac{y}{1-y})]){\mbox e}^{-\theta
\alpha\Scr^{H}_{\alpha}(\frac{y}{1-y})}y^{\alpha\theta-1}{(1-y)}^{-\alpha\theta}
$$
\end{prop}
\Proof  Apply Proposition 3.2 to the density of
$M_{1}(\varrho^{(\alpha\theta)}_{\alpha,H})$ given in Theorem 5.2
\EndProof The next result is also immediate from Theorem 5.2

\begin{prop}Set $\xi/(1-\xi)\overset{d}=M_{1}(H)$, then the
density of
$$A_{\alpha,\xi}\overset{d}=\frac{{[M_{1}(H)]}^{1/\alpha}X_{\alpha}}{{[M_{1}(H)]}^{1/\alpha}X_{\alpha}+1}$$
is
$$
\frac{1}{\pi}\sin(\pi\alpha
[1-\Upsilon_{\alpha}(\frac{y}{1-y})]){\mbox e}^{-
\alpha\Scr^{H}_{\alpha}(\frac{y}{1-y})}y^{\alpha-1}{(1-y)}^{-\alpha-1}
$$\qed
\end{prop}

\Remark We mention briefly the following  example, in the case of
$\alpha=1/2.$ Choose $R^{2}=G_{1}/E,$ then $R^{2}X_{1/2}$ has the
distribution of a Pareto TYPE III law with cdf
$$
\Upsilon_{1/2}(x)=1-{(1+x^{1/2})}^{-1}=\frac{\sqrt{x}}{\sqrt{x}+1}
$$
and ${\mbox e}^{-\theta/2 \Scr^{H}_{1/2}(x)}=\E[{\mbox
e}^{-\frac{1}{x}G_{\theta/2}M_{\theta/2}(\zeta)}]x^{-\theta/2}=x^{\frac{x}{1-x}\frac{\theta}{2}}.$
Where the last expression is obtained from Proposition 5.10 and
Proposition 3.3.\EndRemark

\subsection{Identities for PD local times and skew Bessel Bridges}
We close by showing how the scaling properties discussed in
Proposition 5.7 and 5.8, coupled with the identity
$X_{\alpha\beta}\overset{d}=X_{\alpha}[X_{\beta}]^{1/\alpha},$
translate into interesting identities for skew Bessel bridges and
corresponding local times. The first result is an immediate
consequence of Propositions 4.2 and 5.7.
\begin{prop}For $0<\beta<1$, set $R\overset{d}=X_{\beta},$ that is
$H=\varrho_{\beta},$ then for $\theta>0$, $$ M_{\alpha
\theta}(\varrho_{\alpha\beta})\overset{d}=M_{\alpha\theta}(\varrho_{\alpha}){[M_{\theta}(\varrho_{\beta})]}^{1/\alpha}.
$$ Equivalently,
\begin{enumerate}
\item[(i)]$\ell^{(\alpha\beta,\alpha\theta)}_{S_{\alpha\beta}}\overset{d}=
 \ell^{(\beta,\theta)}_{S_{\beta}}{[\ell^{(\alpha,\alpha\theta)}_{S_{\alpha}}]}^{\beta}$
 \item[(ii)]In terms of the $(\alpha)$-diversity,
 $$
 \frac{S_{\alpha\beta}}{T_{\alpha\beta,\alpha\theta}}\overset{d}=
 \frac{S_{\alpha}}{T_{\alpha,\alpha\theta}}
 {\left[\frac{S_{\beta}}{T_{\beta,\theta}}\right]}^{1/\alpha}
 $$
 which implies that
$T_{\alpha\beta,\alpha\theta}\overset{d}=T_{\alpha,\alpha\theta}{[T_{\beta,\theta}]}^{1/\alpha}.$\qed
 \end{enumerate}
\end{prop}
Now applying Proposition 5.8 leads to,
\begin{prop}For $0<\beta<1$, set $R\overset{d}=X_{\beta},$ that is
$H=\varrho_{\beta},$ which corresponds to
$\varrho_{\alpha,H}=\varrho_{\alpha\beta},$ then for
$0<\alpha\theta<1$, \Eq \label{scaling2}
M_{1}(\varrho^{(\alpha\theta)}_{\alpha\beta})\overset{d}=M_{1}(\varrho^{(\alpha\theta)}_{\alpha}){[M_{\theta}(\varrho_{\beta})]}^{1/\alpha}
\EndEq Recall that
${[\ell^{(\alpha,0)}_{S_{\alpha}}]}^{1/\alpha}\overset{d}=X_{\alpha}\overset{d}=M_{1}(\varrho^{(\alpha)}_{\alpha})\overset{d}=B_{\alpha,1-\alpha}M_{\alpha}(\varrho_{\alpha})$,
then as special cases of~\mref{scaling2}, with the choice of
$\theta=\beta$ and $\theta=1$ respectively,
\begin{enumerate}
\item[(i)]$X_{\alpha \beta}\overset{d}=M_{1}(\varrho^{(\alpha\beta)}_{\alpha\beta})\overset{d}=M_{1}(\varrho^{(\alpha\beta)}_{\alpha}){[M_{\beta}(\varrho_{\beta})]}^{1/\alpha}$
 \item[(ii)]$M_{1}(\varrho^{(\alpha)}_{\alpha\beta})\overset{d}=X_{\alpha}{[M_{1}(\varrho_{\beta})]}^{1/\alpha}.$
\qed
 \end{enumerate}
\end{prop}
The next result shows how to recover $p$-skew Bessel bridges from
$\xi$-skewed bridges.
\begin{prop}For $0<\beta<1$, set $R\overset{d}=c^{\alpha}X_{\beta},$ where $c={(p/q)}^{1/\alpha}$, following Proposition 5.2 this means that $\Lambda_{\alpha,H}=\Lambda_{\alpha\beta,p}.$
Now define $p_{\alpha}=p^{\alpha}/(p^{\alpha}+q^{\alpha}).$ Then,
$\xi_{\theta}\overset{d}=A^{(\beta,\theta)}_{\beta,p_{\alpha}}\overset{d}=M_{\theta}(\Lambda_{\beta,p_{\alpha}})$
and for this choice, \Eq \label{Ascaling}
 A^{(\alpha,\alpha
\theta)}_{\alpha,\xi_{\theta}}\overset{d}=A^{(\alpha\beta,\alpha
\theta)}_{\alpha\beta,p}\overset{d}=M_{\alpha\theta}(\Lambda_{\alpha\beta,p})\overset{d}=P_{\alpha\beta,\alpha\theta}(C)
\EndEq Results recovering $p$-skew  Bessel bridges follows as
special cases of~\mref{Ascaling}, with the choice of
$\theta=\beta$ and $\theta=1$ respectively,
\begin{enumerate}
\item[(i)]$A^{(\alpha,\alpha
\beta)}_{\alpha,\xi_{\beta}}\overset{d}=A^{(\alpha\beta,\alpha
\beta)}_{\alpha\beta,p}$
\item[(ii)]$A^{(\alpha,\alpha
)}_{\alpha,\xi_{1}}\overset{d}=A^{(\alpha\beta,\alpha
)}_{\alpha\beta,p}$
\end{enumerate}
\end{prop}
\Proof We present a proof for clarity. First, as noted, the choice
of $R=c^{\alpha}X_{\beta}$ coupled with proposition 5.2 shows that
$\Lambda_{\alpha,H}=\Lambda_{\alpha\beta,p},$ leading to
~\mref{Ascaling}. What remains is to verify that
$\xi_{\theta}\overset{d}=A^{(\beta,\theta)}_{\beta,p_{\alpha}}.$
Noting Proposition 5.4,
$\xi_{\theta}\overset{d}=M_{\theta}(H^{(c^{\alpha},c^{\alpha})})$
, but here
$H^{(c^{\alpha},c^{\alpha})}=\Lambda_{\beta,p_{\alpha}},$ when
$H=\varrho_{\beta}.$ \EndProof

So in closing we see, as a special case, that occupation time of a
Brownian bridge up to time $1$, say
$A^{br}_{1/2}:=A^{(1/2,1/2)}_{1/2,1/2}\overset{d}=B_{1/2,1/2}$ is
equivalent in distribution to the time time spent positive up to
time $1$ of a $A^{(\beta,\beta)}_{\beta,1/2}$ randomly skewed
process $B^{(\alpha,1/2)}_{\alpha}(t)$ for any $\alpha=1/(2\beta)$
and $1/2<\beta<1.$ That is to say, this process is randomly skewed
by the random variable corresponding to the time spent positive up
till time $1$ of a Bessel bridge of dimension $2-2\beta<1,$
$B^{(\beta,\beta)}_{\beta}(t).$

\subsection*{Acknowledgement}I would like to thank Marc Yor for
stimulating conversations and exchanges of ideas related to this
work. I also wish to thank Jessica Tressou for help in clarifying
some ideas.

\vskip0.2in \centerline{\Heading References} \vskip0.2in

\tenrm
\def\smc{\tensmc}
\def\sl{\tensl}
\def\bf{\tenbold}
\baselineskip0.15in

\Ref \by Aldous, D., and Pitman, J.\yr 2004 \paper Two recursive
decompositions of Brownian bridge related to the asymptotics of
random mappings. In \emph{In Memoriam Paul-Andre Meyer -
S\'eminaire de Probabilit\'es XXXIX.} (Yor, M. and \'Emery, M.,
Eds.), 269-303, Springer Lecture Notes in Math. \textbf{1874}.
Springer, Berlin. \EndRef

\Ref  \by Barlow, M., Pitman, J. and Yor, M. \yr 1989 \paper Une
extension multidimensionnelle de la loi de l'arc sinus. In
\textit{S\'eminaire de Probabilit\'es XXIII} (Azema, J., Meyer,
P.-A. and Yor, M., Eds.), 294--314, Lecture Notes in Mathematics
\textbf{1372}. Springer, Berlin.\EndRef

\Ref \by Bertoin, J. \yr 2006 \book Random fragmentation and
coagulation processes \publ Cambridge University Press\EndRef

\Ref \by Bertoin, J., Fujita, T., Roynette, B., and Yor, M. \yr
2006 \paper On a particular class of self-decomposable random
variables: the duration of a Bessel excursion straddling an
independent exponential time. To appear in Prob. Math. Stat
\EndRef

\Ref \by Bertoin, J. and Yor, M. \yr 1996 \paper  Some
independence results related to the arc-sine law. \emph{J.
Theoret. Probab.} \textbf{9} 447-458\EndRef

\Ref \by Bondesson, L. \yr 1992 \paper Generalized gamma
convolutions and related classes of distributions and densities.
Lecture Notes in Statistics, 76. Springer-Verlag, New York \EndRef

\Ref \by Chaumont, L. and Yor, M. \yr 2003 \book Exercises in
probability. A guided tour from measure theory to random
processes, via conditioning. Cambridge Series in Statistical and
Probabilistic Mathematics, 13 \publ Cambridge University Press,
Cambridge\EndRef

 \Ref \by    Cifarelli, D. M. and Melilli, E. \yr    2000
\paper Some new results for Dirichlet priors \jour \AnnStat \vol
28 \pages 1390-1413\EndRef

\Ref \by Cifarelli, D. M. and Regazzini, E.\yr 1990 \paper
Distribution functions of means of a Dirichlet process. {\it Ann.
Statist.} \textbf{18}, 429--442 (Correction in \textit{Ann.
Statist.} (1994) \textbf{22}, 1633-1634)\EndRef

\Ref \by Devroye, L. \yr 1990 \paper A note on Linnik's
distribution \jour Statist. Probab. Lett. \vol 9 \pages
305-306\EndRef

\Ref \by Devroye, L. \yr 1996 \paper Random variate generation in
one line of code, in: 1996 Winter Simulation Conference
Proceedings, ed. J.M. Charnes, D.J. Morrice, D.T. Brunner and J.J.
Swain, pp. 265-272, ACM\EndRef

\Ref \by  Diaconis, P. and Kemperman, J. \yr    1994 \paper Some
new tools for Dirichlet priors. Bayesian Statistics 5 (J.M.
Bernardo, J.O. Berger, A.P. Dawid and A.F.M. Smith eds.), Oxford
University Press, pp. 97-106 \EndRef

\Ref \by    Ferguson, T. S. \yr    1973 \paper A Bayesian analysis
of some nonparametric problems \jour  \AnnStat \vol   1 \pages
209-230 \EndRef

\Ref\by Fujita, T. and Yor, M. \yr 2006 \paper An interpretation
of the results of the BFRY paper in terms of certain means of
Dirichlet Processes. preprint\EndRef

\Ref \by Hjort, N.L. and Ongaro, A. \yr 2005 \paper Exact
inference for random Dirichlet means. \emph{Stat. Inference Stoch.
Process.}, \textbf{8} 227-254 \EndRef

\Ref \yr  2001 \by Ishwaran,~H. and James,~L.~F. \paper Gibbs
sampling methods for stick-breaking priors \jour Journal of the
American Statistical Association \vol 96 \pages 161-173 \EndRef

\Ref \yr    2005 \by James,~L.~F. \paper Functionals of Dirichlet
processes, the Cifarelli-Regazzini identity and Beta-Gamma
processes \jour Annals of Statistics \vol 33 \pages 647-660\EndRef

\Ref \yr 2006 \by James, L.F. \paper Laws and likelihoods for
Ornstein Uhlenbeck-Gamma and other BNS OU stochastic volatilty
models with extensions. http://arxiv.org/abs/math/0604086 \EndRef

\Ref \yr 2006 \by James, L.F., Lijoi, A. and I. Pr\"unster \paper
Distributions of functionals of the two parameter
Poisson-Dirichlet process.
http://arxiv.org/abs/math.PR/0609488\EndRef

\Ref \by Kasahara, Y. and  Watanabe, S.\yr 2005 \paper Occupation
time theorems for a class of one-dimensional diffusion processes.
\emph{Period. Math. Hungar.} \textbf{50}, 175-188\EndRef

 \Ref \by Lamperti, J. \yr 1958 \paper An occupation time
theorem for a class of stochastic processes \jour Trans. Amer.
Math. Soc.\vol 88 \pages 380-387\EndRef

\Ref \by L\'evy, P. \yr 1939 \paper Sur certains processus
stochastiques homog\'enes \jour Compositio Math. \vol 7 \pages
283-339\EndRef

\Ref \by Perman, M., Pitman, J. and Yor, M. \yr 1992 \paper
Size-biased sampling of Poisson point processes and excursions.
\emph{Probab. Theory Related Fields.} \textbf{92}, 21-39\EndRef

\Ref \by Pitman, J. \yr 1996 \paper Some developments of the
Blackwell-MacQueen urn scheme. In {\it Statistics, Probability and
Game Theory} (T.S. Ferguson, L.S. Shapley and J.B. MacQueen, Eds.)
245--267. IMS Lecture Notes-Monograph series\EndRef

\Ref \by Pitman, J. \yr 1999 \paper Coalescents with multiple
collisions. \jour \AnnProb \vol 27 \pages 1870-1902\EndRef

\Ref \by Pitman, J. \yr 2006 \paper \emph{Combinatorial Stochastic
Processes.} Ecole d'Et\'e de Probabilit\'es de Saint-Flour XXXII –
2002. Lecture Notes in Mathematics \textbf{1875}. Springer, Berlin
\EndRef

\Ref  \by Pitman, J., and  Yor, M. \yr 1992 \paper Arcsine laws
and interval partitions derived from a stable subordinator.
\emph{Proc. London Math. Soc.} \textbf{65} 326-356 \EndRef

\Ref  \by Pitman, J., and  Yor, M. \yr 1997a \paper The
two-parameter Poisson-Dirichlet distribution derived from a stable
subordinator \jour \AnnProb \vol 25 \pages 855-900\EndRef

\Ref \by Pitman, J., and  Yor, M. \yr 1997b \paper On the relative
lengths of excursions derived from a stable subordinator. In:
\textit{S\'eminaire de Probabilit\'es XXXI} (Azema, J., Emery, M.
and Yor, M., Eds.), 287--305, Lecture Notes in Mathematics
\textbf{1655}. Springer, Berlin\EndRef

\Ref \by Vershik, A., Yor, M. and Tsilevich, N. \yr 2001 \paper On
the Markov-Krein identity and quasi-invariance of the gamma
process. \emph{Zap. Nauchn. Sem. S.-Peterburg. Otdel. Mat. Inst.
Steklov. }(POMI) \textbf{283} 21--36. [In Russian. English
translation in \emph{J. Math. Sci.} \textbf{121} (2004),
2303--2310]\EndRef

\Ref \by Watanabe, S. \yr 1995 \paper Generalized arc-sine laws
for one-dimensional diffusion processes and random walks.
Stochastic analysis (Ithaca, NY, 1993), 157-172, Proc. Sympos.
Pure Math., 57, Amer. Math. Soc., Providence, RI \EndRef

\Tabular{ll}

Lancelot F. James\\
The Hong Kong University of Science and Technology\\
Department of Information and Systems Management\\
Clear Water Bay, Kowloon\\
Hong Kong\\
\rm lancelot\at ust.hk\\
\EndTabular

\end{document}